\documentclass[11pt]{article}
\usepackage[utf8]{inputenc}
\usepackage{amsmath,amssymb,epsfig,bbm}
\usepackage{stmaryrd,mathabx}
\usepackage{comment}
\usepackage{color}
\usepackage[T1]{fontenc}

\usepackage{interval}

\usepackage[textsize=small]{todonotes}
\usepackage{enumitem}
\usepackage{varwidth}
\setlist{nolistsep}
\usepackage{hyperref}

\newtheorem{defi}{Definition}[section]
\newtheorem{prop}[defi]{Proposition}
\newtheorem{theo}[defi]{Theorem}
\newtheorem{theofr}[defi]{Théorème}
\newtheorem{conj}[defi]{Conjecture}
\newtheorem{lemm}[defi]{Lemma}
\newtheorem{lemmfr}[defi]{Lemme}
\newtheorem{coro}[defi]{Corollary}
\newtheorem{rema}[defi]{Remark}
\newtheorem{exem}[defi]{Example}
\newtheorem{exems}[defi]{Examples}

\newcommand{\bdefi}{\begin{defi}}
\newcommand{\edefi}{\end{defi}}
\newcommand{\bprop}{\begin{prop}}
\newcommand{\eprop}{\end{prop}}
\newcommand{\btheo}{\begin{theo}}
\newcommand{\etheo}{\end{theo}}
\newcommand{\btheofr}{\begin{theofr}}
\newcommand{\etheofr}{\end{theofr}}
\newcommand{\blemm}{\begin{lemm}}
\newcommand{\elemm}{\end{lemm}}
\newcommand{\blemmfr}{\begin{lemmfr}}
\newcommand{\elemmfr}{\end{lemmfr}}
\newcommand{\brema}{\begin{rema}}
\newcommand{\erema}{\end{rema}}
\newcommand{\bexer}{\begin{exem}}
\newcommand{\eexer}{\end{exem}}
\newcommand{\bexems}{\begin{exems}}
\newcommand{\eexems}{\end{exems}}
\newcommand{\bconj}{\begin{conj}}
\newcommand{\econj}{\end{conj}}
\newcommand{\bcoro}{\begin{coro}}
\newcommand{\ecoro}{\end{coro}}
\newcommand{\dem}{\noindent{\bf Proof. }}

\usepackage{mathrsfs}
\renewcommand\mathcal{\mathscr}

\newcommand{\D}{{\cal D}}

\newcommand{\F}{{\cal F}}

\renewcommand{\H}{{\cal H}}

\renewcommand{\L}{{\cal L}}

\newcommand{\R}{{\cal R}}
\newcommand{\maths}[1]{{\mathbb #1}}  

\newcommand{\CC}{\maths{C}}

\newcommand{\HH}{\maths{H}}

\newcommand{\NN}{\maths{N}}

\newcommand{\QQ}{\maths{Q}}
\newcommand{\RR}{\maths{R}}

\newcommand{\ZZ}{\maths{Z}}

\newcommand{\weakstar}{\overset{*}\rightharpoonup}
\newcommand{\ra}{\rightarrow}
\newcommand{\bs}{\backslash}

\newcommand{\wt}[1]{{\widetilde{#1}}}

\newcommand{\ga}{\gamma}
\newcommand{\Ga}{\Gamma}

\newcommand{\cqfd}{\hfill$\Box$}

\newcommand{\bigO}{\operatorname{O}}
\newcommand{\BV}{\operatorname{BV}}
\newcommand{\card}{{\operatorname{Card}}}
\newcommand{\CAT}{\operatorname{CAT}}

\newcommand{\id}{\operatorname{id}}
\renewcommand{\Im}{{\operatorname{Im}}}

\newcommand{\Leb}{\operatorname{Leb}}
\renewcommand{\log}{\operatorname{ln}}

\newcommand{\mult}{\operatorname{mult}}
\newcommand{\ols}{\operatorname{OL}}

\newcommand\Perp{\operatorname{Perp}}

\renewcommand{\Re}{{\operatorname{Re}}}

\newcommand{\sg}{\operatorname{sg}}

\newcommand{\Var}{\operatorname{Var}}

\newcommand{\hdr}{{\HH}^2_\RR}

\newcommand{\PSL}{\operatorname{PSL}}

\newcommand{\PSLZ}{\operatorname{PSL}_{2}(\ZZ)}

\newcommand{\n}{\operatorname{\tt n}}

\newcounter{fig}

\title{Pair correlations of logarithms of integers}
\author{Jouni Parkkonen \and Fr\'ed\'eric Paulin} 
\date{With an appendix by \'Etienne Fouvry\\\mbox{}\\ \today}

\begin{document}
\bibliographystyle{../alphanum}
\maketitle
\begin{abstract} 
We study the correlations of pairs of logarithms of positive integers
at various scalings, either with trivial weigths or with weights given
by the Euler function, proving the existence of pair correlation
functions. We prove that at the linear scaling, the pair correlations
exhibit level repulsion, as it sometimes occurs in statistical
physics. We prove total loss of mass phenomena at superlinear
scalings, and Poissonian behaviour at sublinear scalings. The case of
Euler weights has applications to the pair correlation of the lengths
of common perpendicular geodesic arcs from the maximal Margulis cusp
neighborhood to itself in the modular curve $\PSL_2(\ZZ)\bs\hdr$.
\end{abstract}

\section{Introduction}
\label{sec:intro}

When studying the asymptotic distribution of a sequence of finite
subsets of $\RR$, finer information is sometimes given by the
statistics of the spacing (or gaps) between pairs or $k$-tuples of
elements, seen at an appropriate scaling. This problematic is largely
developped in quantum chaos, including energy level spacings or
clusterings, and in statistical physics, including molecular repulsion
or interstitial distribution.  The general setting of such a study may
be described as follows.  Let $\F=(\,(F_N)_{N\in\NN}, \;\omega)$ be a
nondecreasing sequence of finite subsets $F_N$ of a finite dimensional
Euclidean space $E$, endowed with a {\it multiplicity function}
$\omega: \bigcup_{N\in\NN} F_N \ra\; ]\,0,+\infty\,[\,$ (or {\it
weight} function). Let $\psi:\NN\ra\; ]\,0,+\infty\,[$ be a
nondecreasing {\it scaling function}. The {\it pair correlation
measure of $\F$ at time $N$ with scaling $\psi(N)$} is the measure
on $E$ with finite support
$$
\R^{\F,\psi}_N=\sum_{x,y\in F_N\;:\;x\neq y}\;\omega(x)\,\omega(y)\,
\Delta_{\psi(N)(y-x)}\;,
$$
where $\Delta_z$ denotes the unit Dirac mass at $z$.  When the
sequence of measures $(\R^{\F,\psi}_N)_{N\in\NN}$, appropriately
renormalized, weak-star converges to a measure $g\,\Leb_E$ absolutely
continuous with respect to the Lebesgue measure $\Leb_E$ of $E$, the
Radon-Nikodym derivative $g=g_{\F,\psi}$ is called the asymptotic {\it
  pair correlation function} of $\F$ {\em for the scaling}
$\psi$. When $g_{\F,\psi}$ vanishes on a neighbourhood of $0$ in $E$,
we say that the pair $(\F,\psi)$ {\it exhibits level repulsion}.

If the family $\F$ consists of subsets of the unit interval $[0,1]$,
then it is customary to use the cardinality of the finite set $F_N$ as
the scaling function. See for example \cite{BocZah05}, where
$F_N=\{\frac{p}{q}:p,q\in\NN, p\leq q, (p,q)=1, 0<q\leq N\}$ is the
set of Farey fractions of order $N$ in $[0,1]$ (without
multiplicities, hence $\omega\equiv 1$), so that $\psi(N)=
\frac{3N^2}{\pi} +\bigO(N\ln N)$.  Montgomery studied (under the
Riemann hypothesis) the pair correlations of the imaginary parts of
the zeros (with their multiplicity as zeros) of the Riemann zeta
function $\zeta$ in the seminal paper \cite{Montgomery73}. The number
of zeros $\frac 12+i t$ of $\zeta$ with imaginary part $t$ in the
interval $\interval 0N$ is asymptotic to $\frac{N\ln N}{2\pi}$ as
$N\to+\infty$ and the scaling used in \cite{Montgomery73} is,
analogously to the unit interval case of \cite{BocZah05}, by the
inverse of the {\em average gap} : $\psi(N)=(\frac{N\ln
N}{2\pi})/N=\frac{\ln N}{2\pi}$.
%See also \cite{BleShiZel00CMP,BleShiZel00Inv}.

\medskip
In Sections \ref{sect:logpaircorrel} and
\ref{sect:logpaircorrelscaled}, we study the pair correlations of the
family of the logarithms of positive integers
$$
\L_\NN=\big(\,(L_N=\{\ln n\;:\;0<n\leq N\})_{N\in\NN},\omega\equiv 1\big)
$$
without multiplicities. In order to simplify the statements in this
introduction, we only consider power scalings $\psi:N\mapsto N^\alpha$ for
$\alpha\ge 0$, and we denote these scaling functions by $\id^\alpha$.

\btheo\label{theo:intro1} Let $\alpha\ge 0$.  As $N\ra+\infty$, the
normalized pair correlation measures $\frac{1}{N^{2-\alpha}}\;
\R_N^{\,\L_\NN,\,\id^\alpha}$ on $\RR$ weak-star converge to a measure
$g_{\L_\NN,\,\id^\alpha}\;\Leb_\RR$ with pair correlation function given by
$$
g_{\L_\NN,\,\id^\alpha}:t\mapsto\begin{cases}
\frac{1}{2}\;e^{-|s|} & \textrm{if } \alpha=0\\
\frac{1}{2}&\textrm{if } 0<\alpha<1\\
\frac{1}{2\,t^2}\;\lfloor |t|\rfloor
\big(\,\lfloor |t|\rfloor+1\big)
     &\textrm{if }\alpha=1\\
     0&\textrm{if } \alpha>1\;.
\end{cases}
$$
\etheo

We refer to Theorems \ref{theo:logpaircorrel} and
\ref{theo:logpaircorr}, for more complete versions of Theorem
\ref{theo:intro1}, with congruence restrictions and with more general
scaling functions, as well as for error terms. These error terms, as
well as the ones in Theorems \ref{theo:logpaircorrelphi} and
\ref{theo:logpaircorrelphipsi}, constitute the main technical parts of
this paper.

The renormalisation by $\frac{1}{N^{2-\alpha}}$ in Theorem
\ref{theo:intro1} is naturally chosen in order for the pair
correlation function to be finite.  As the finite set $L_N$, whose
order is $N$, is contained in the minimal interval $\interval 0{\ln
  N}$, the average gap in $L_N$ is $\frac{\log N}N$. Scaling by the
inverse $\psi(N)= \frac{N}{\ln N}$ of the average gap (as in a particular
case of Theorem \ref{theo:logpaircorr}), as well as by $\id^\alpha$
for $0<\alpha<1$ (as in the above statement) gives a pair correlation
function which is constant nonzero, a characteristic of a {\em
  Poissonian distribution}. As in the above result for $\alpha>1$ and
more generally by Theorem \ref{theo:logpaircorr}, if the scaling
function $\psi$ grows faster than linearly, then the pair correlation
function vanishes : the empirical measures $\R_N^{\,\L_\NN,\,\psi}$
have a total loss of mass at infinity, actually whatever the
renormalisation is (the support of the measure itself converges to
infinity).  The transition from Poissonian to zero correlation occurs
at linear scalings, where a more exotic pair correlation function
appears (see for instance \cite{RudSar98} for Poissonian and
\cite{HofKal21,LarSto20} for non-Poissonian pair correlation
phenomena).  Since $g_{\L_\NN,\,\id}$ vanishes on $]-1,1[$, the pair
    $(\L_\NN,\,\id)$ exhibits a level repulsion.

The figure below gives the graph of the pair correlation function
$g_{\L_\NN,\,\id}$ of $\L_\NN$ at the linear scaling $\psi=\id:
N\mapsto N$ in the interval $[-15,15]$ compared with the graph of the
constant function $\frac 12$. The graph is similar to certain radial
distribution functions in statistical physics, see for example
\cite[Sect.~II]{ZerPri27}, \cite[Fig.~7]{SanLop16}, \cite[page
  199]{Chandler87} or \cite[page 18]{Bohigas91}.

\begin{center}
\includegraphics[width=13cm]{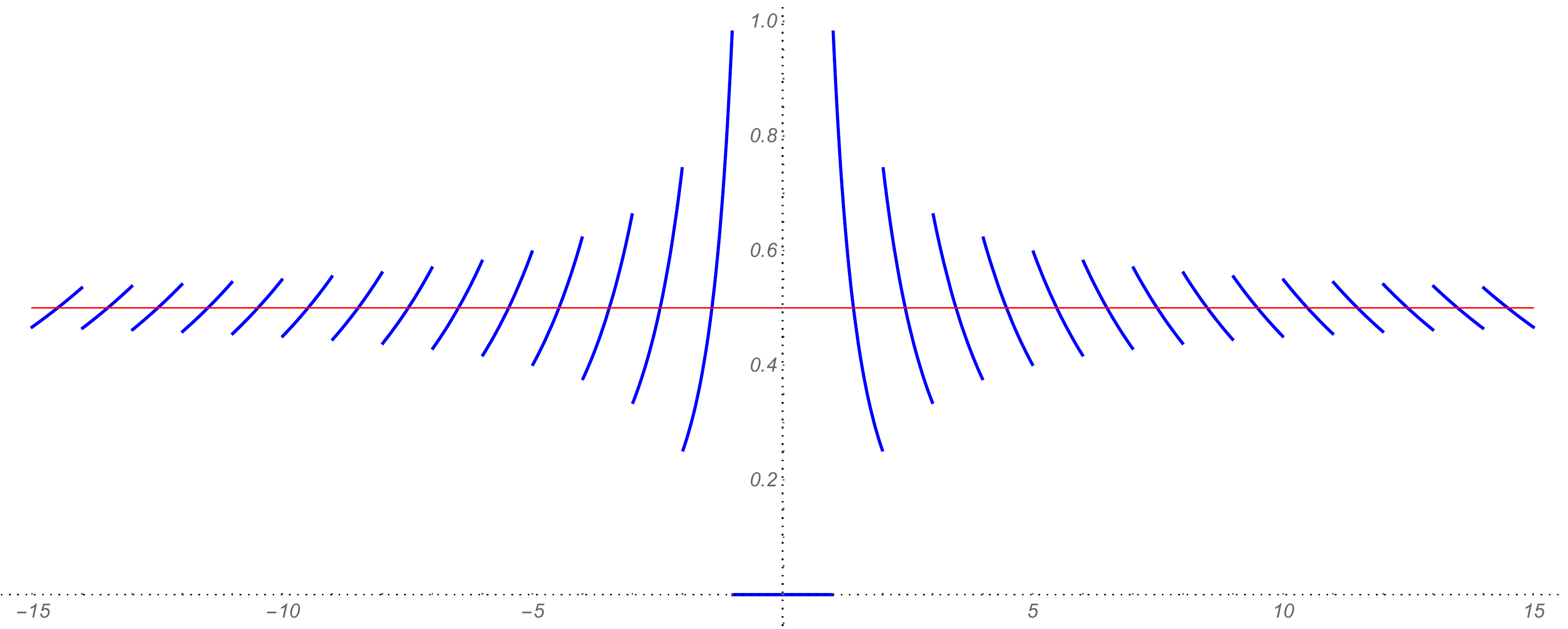}
\end{center}

Instead of the pair correlations, one can study the gaps between
consecutive elements in the subsets $F_N$ of the real line or, most
often, of the unit interval.  Marklof and Str\"ombergsson
\cite{MarStr13} have computed the gap distribution of the fractional
parts of the family $\L_\NN$ (with a linear scaling and linear
renormalisation) and showed that the limiting gap distribution has two
jump discontinuities.

\medskip
In Section \ref{sect:geometricmotivation}, we prove that the pair
correlation measures of the lengths of the common perpendiculars
between the maximal Margulis cusp neighbourhood and itself in the
modular curve $\PSLZ\bs\hdr$ are (up to a factor $2$) the pair
correlation measures of the weighted family
$$
\L_\NN^\varphi=\big(\,(L_N=\{\ln n\;:\; 0<n\leq N\})_{N\in\NN},\;
\omega=\varphi\circ \exp\big)
$$
of logarithms of integers, with weights given by the Euler
function $\varphi:n\mapsto \card(\ZZ/n\ZZ)^\times$, see Proposition
\ref{prop:corrhoromodular}. See \cite{PolSha06,PolSha13} for results
on the pair correlation of the lengths of closed geodesics in
negatively curved manifolds.

We study the pair correlations of the arithmetically defined family
$\L_\NN^\varphi$ in Sections \ref{sect:weightlogpaircorrel} and
\ref{sect:weightscaledlogpaircorrel}, where we find the pair
correlation function without scaling and with linear scaling.

\btheo\label{theo:intro2} (1) As $N\ra+\infty$, the pair correlation
measures $\R^{\L_\NN^\varphi,1}_N$ on $\RR$, renormalized to be
probability measures, weak-star converge to the probability measure
$g_{\L_\NN^\varphi,1}\; \Leb_\RR$, with pair correlation function
$g_{\L_\NN^\varphi,1}:s\mapsto e^{-\,2\,|s|}$.

\smallskip\noindent (2) As $N\ra+\infty$, the normalized pair
correlation measures $\frac{1}{N^3}\,\R^{\L_\NN^\varphi,\,\id}_N$
(with linear scaling) on $\RR$ weak-star converge to the measure
$g_{\L_\NN^\varphi,\,\id}\; \Leb_\RR$, with pair correlation function
\begin{equation}\label{eq:gLvarphiid}
g_{\L_\NN^\varphi,\,\id}: s\mapsto  \frac{1}{s^4}
\prod_{p {\rm ~prime}}(1-\frac{2}{p^2})\;\sum_{k=1}^{\lfloor\,|s|\,\rfloor} \;k^3
\prod_{p {\rm ~prime},\; p\,\mid\, k}(1+\frac{1}{p(p^2-2)})\;.
\end{equation}
\etheo

We refer to Theorems \ref{theo:logpaircorrelphi} and
\ref{theo:logpaircorrelphipsi} for more complete versions of Theorem
\ref{theo:intro2} with congruence restrictions, and for error
terms. When the congruences are nontrivial, the proof of the second
claim of Theorem \ref{theo:intro2} uses a generalization of Mirsky's
formula (see \cite{Mirsky49}) that is proved in Appendix
\ref{appendixFouvry} by \'Etienne Fouvry.

The figure below gives the graph of the pair correlation function
$g_{\L_\NN^\varphi,\,\id}$ compared with the graph of the constant
function with value $\frac{1}{4}{\displaystyle\prod_{p {\rm ~prime}}}
\big(1-\frac{2}{p^2}\big) \big(1+\frac{1}{p^2(p^2-2)}\big)\simeq
0.09239$, which is the limit of the pair correlation function
$g_{\L_\NN^\varphi,\,\id}$ at $\pm\infty$ by Proposition
\ref{prop:Fouvry2}.

\begin{center}
\includegraphics[width=13cm]{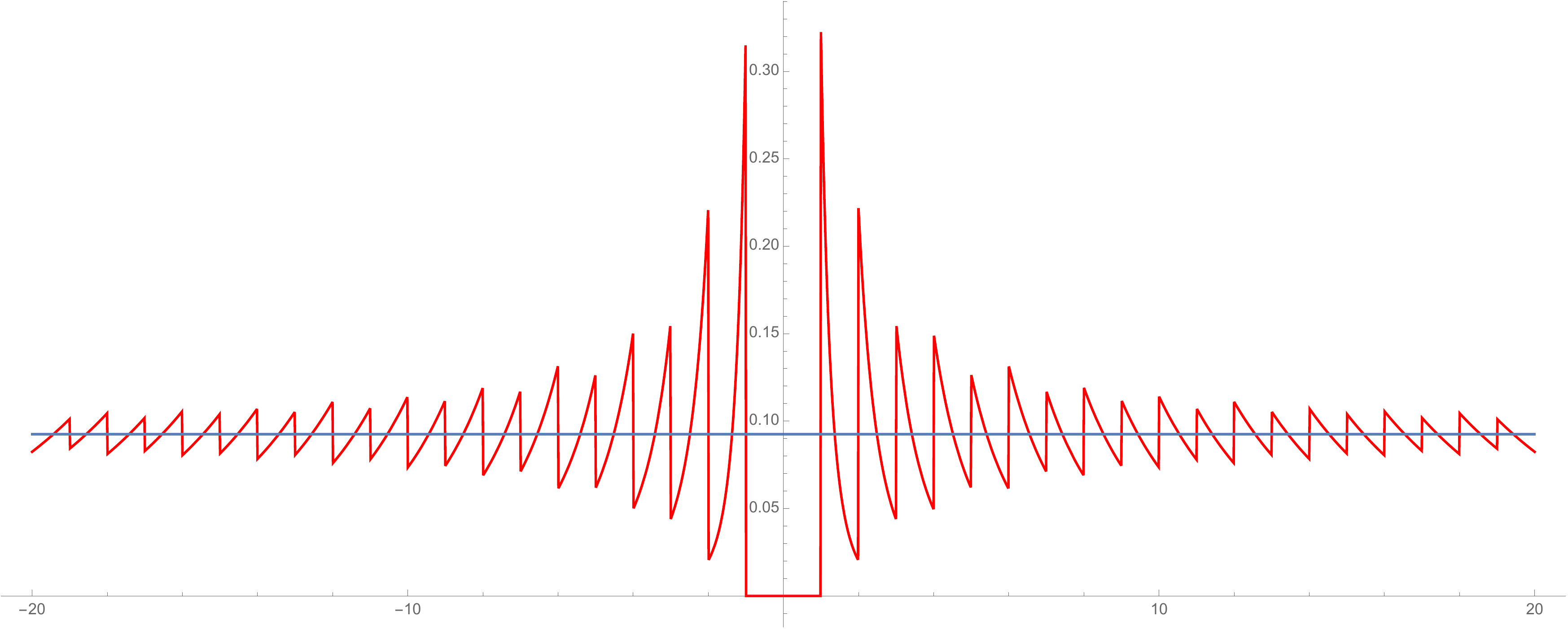}
\end{center}

Theorems \ref{theo:logpaircorrelphi} and
\ref{theo:logpaircorrelphipsi} imply pair correlation results for the
common perpendiculars of cusps neighborhoods in the modular curve and
on quotients of the hyperbolic plane by Hecke congruence subgroups of
$\PSL(\ZZ)$, see Corollary \ref{coro:comperphoromodular} for precise
statements.

\medskip\noindent {\bf Further directions.}
It would be interesting, given a discrete subgroup $\Ga$ of
$\PSL_2(\RR)$, to study the asymptotic of the pair correlation
measures of the complex translation lengths $\ell_\CC(\ga)$ with
absolute value at most $N$ of the elements $\ga\in\Ga$, and given a
discrete subgroup $\Ga$ of a semi-simple connected real Lie group $G$
with finite center and without compact factor, of the Cartan
projections $\mu(\ga)$ with Killing norm at most $N$ of the elements
$\ga\in\Ga$.  See Section \ref{sect:geometricmotivation} for the
problem of the asymptotic of the pair correlation measures of common
perpendiculars in negative curvature, which will be studied more
completely in subsequent works of the authors.

When the finite-dimensional Euclidean space $E$ (where the family
of finite sets $(F_N)_{N\in\NN}$ sits) is replaced by a locally compact
metric space $(X,d)$, we may also consider the positive measure on
$]0,+\infty[$ with finite support $\R^{\F,\psi}_N= \sum_{x,y\in
F_N\;:\;x\neq y}\;\omega(x)\,\omega(y)\,\Delta_{\psi(N)d(x,y)}$.

\medskip
\noindent{\small {\it Acknowledgements: } The authors thank a lot
  Etienne Fouvry for his proofs of Lemma \ref{lem:Fouvry1},
  Proposition \ref{prop:Fouvry2} and Theorem \ref{theo:mirskyfouvry}
  and for the agreeing to contribute the appendix to this paper.  This
  research was supported by the French-Finnish CNRS IEA BARP. }
  
\medskip\noindent {\bf Notation.} We introduce here some of the
notation used throughout the paper.
  
The pushforward of a measure $\mu$ by a mapping $f$ is denoted by
$f_*\mu$. We denote by $\sg:\RR\ra\RR$ the change of sign map
$t\mapsto -t$.
  
For every interval $I$ in $\RR$, we denote by $\Leb_I$ the Lebesgue
measure on $I$ and by $\mathbbm{1}_I$ the characteristic function of
$I$. We denote by $\BV(I)$ the vector space of measurable functions
$f:I\ra\RR$ with finite total variation $\Var(f)$.  For every
$k\in\{0,1\}$, we denote by $C^k_{\rm c}(I)$ the real vector space of
$C^k$-smooth functions $f:I\ra \CC$ with compact support in $I$.  We
denote by $\|f\|_\infty =\sup_{x\in I} |f(x)|$ the uniform norm of
$f\in C^0_{\rm c}(I)$.

In addition to the above, more or less standard, notation, we will
use the following indexing sets in Sections \ref{sect:logpaircorrel},
\ref{sect:logpaircorrelscaled}, \ref{sect:weightlogpaircorrel} and
\ref{sect:weightscaledlogpaircorrel}.  Let us fix throughout the paper
$a,b\in\NN-\{0\}$ with $a\leq b$.  For every $N\in\NN-\{0\}$, let
\begin{align*}
  I_N=I_{N,a,b}&=\{(m,n)\in\NN^2\;:\;
  0<m,n\leq N,\;m\neq n,\;\;m,n\equiv a\!\!\mod b\}\;,\\
  I^-_N&=\{(m,n)\in\NN^2\;:\; 0<m<n\leq N,\;\;m,n\equiv a\!\!\mod b\}\,\\
  I^+_N&=\{(m,n)\in\NN^2\;:\; 0<n<m\leq N,\;\;m,n\equiv a\!\!\mod b\}\;,
\end{align*}
so that $I_N=I^-_N\sqcup I^+_N$ is the disjoint union of $I^-_N$ and
$I^+_N$.

Except in the proof of Lemma \ref{lem:Fouvry1}, of Proposition
\ref{prop:Fouvry2}, and in the whole Appendix \ref{appendixFouvry},
for every function $g$ of a variable in $\NN-\{0\}$, possibly
depending on parameters (including $a$ and $b$), we will denote by
$\bigO(g)$ (or $\bigO_b(g)$ when we want to insist on the possible
dependence on the parameter $b$) any function $f$ on $\NN-\{0\}$ such
that there exists a constant $C'$ depending only on the parameter $b$
and a constant $N_0$ possibly depending on the parameters such that
for every $N\geq N_0$, we have $|f(N)|\leq C\;|g(N)|$.

\section{Pair correlations without weights nor scaling}
\label{sect:logpaircorrel}

For every $N\in\NN-\{0\}$, the (not normalised) {\em pair correlation
  measure} of the logarithms of integers congruent to $a$ modulo $b$
at time $N$, with trivial multiplicities and with trivial scaling
function, is
$$
\nu_N =\sum_{(m,\,n)\in I_N}\;\Delta_{\ln m-\ln n}\;.
$$
If we consider the following nondecreasing sequence of finite
subsets of $\RR$ with trivial multiplicity
$$
\L_\NN^{a,b}=\big(\,(L_N^{a,b}=\{\ln n : 0<n\leq N,\;\;
n\equiv a\!\!\mod b\})_{N\in\NN},\;\omega\equiv 1\big)\;,
$$
then, with the notation of the introduction, we have
$\L_\NN^{1,1}=\L_\NN$ and $\nu_N=\R^{\L_\NN^{a,b},1}_N$.

\btheo\label{theo:logpaircorrel} As $N\ra+\infty$, the measures
$\nu_N$ on $\RR$, renormalized to be probability measures, weak-star
converge to the measure absolutely continuous with respect to the
Lebesgue measure on $\RR$, with Radon-Nikodym derivative the function 
$g_{\L_\NN^{a,b},1}:s\mapsto\frac{1}{2}\;e^{-|s|}$:
$$
\frac{\nu_N}{\|\nu_N\|}\;\;\;\weakstar\;\;\;g_{\L_\NN^{a,b},1}\;\Leb_\RR\;.
$$
Furthermore, for every $f\in C^0_{\rm c}(\RR)\cap\BV(\RR)$, we have
$$
\frac{\nu_N}{\|\nu_N\|}(f)=\frac{1}{2}\;\int_{s\in\RR}\;f(s)\;e^{-|s|}\;ds
\;+\;\bigO_b\Big(\frac{\|f\|_\infty+\Var(f)}{N}\Big)\;.
$$
\etheo

When $a=b=1$, this result implies the case $\alpha=0$ of Theorem
\ref{theo:intro1} in the introduction, with pair correlation function
$g_{\L_\NN,1}=g_{\L_\NN^{1,1},1}$.

\medskip
\dem For every $q\in\NN-\{0\}$ with $q\equiv a\!\!\mod b$, let
$q'\in\NN$ be such that $q=a+q'b$ and
\begin{equation}\label{eq:defJsubq}
J_q=\{p\in\NN :  0<p<q,\;p\equiv a\!\!\mod b\}=
\{a+kb : 0\leq k< q'\}\;.
\end{equation}
%so that $I_n=\coprod_{0<q\leq N\;:\:q\equiv a\!\!\mod  b}
%\{(p,q)\in\NN^2 : p\in J_q\}$.
Let
$$
\omega_q=\sum_{p\in J_q}\Delta_{\frac{p}{q}}\;,
$$
which is a finitely supported measure on $[0,1]$, with total mass
$\|\omega_q\| =q'$. When $q'\neq 0$, we hence have $\|\omega_q\|
=\frac{q}{b}+\bigO(1)$ and $\frac{1}{\|\omega_q\|} =\frac{b}{q}
+\bigO(\frac{1}{q^2})$. When $q'\neq 0$, we denote by
$\overline{\omega_q} =\frac{\omega_q}{\|\omega_q\|}$ the
renormalisation of $\omega_q$ to a probability measure on $[0,1]$. By
well known Riemann sum arguments, we have, as $q\ra+\infty$,
$$
\overline{\omega_q}\;\;\weakstar\;\;\Leb_{[0,1]}\;.
$$
Let $f\in \BV([0,1])$, and note that $f$ is bounded, with
$\|f\|_\infty \leq f(0)+\Var(f)$. Denoting by $M_k$ and $m_k$ the
maximum and minimum respectively of $f$ on $[\frac{a+kb}{q},
  \frac{a+(k+1)b}{q}]$ for $0\leq k< q'$, we have
\begin{align*}
  &\Big|\int_0^1f(t)\,dt-\frac{b}{q}\,\omega_q(f)\,\Big|=
  \Big|\int_0^1f(t)\,dt-\sum_{p\in J_q}\;\frac{b}{q}\;f(\frac{p}{q})
  \;\Big|\\ \leq\;&  \Big|\int_0^{\frac{a}{q}}f(t)\,dt\;\Big|
  +\sum_{k=0}^{q'-1}\Big|\int_{\frac{a+kb}{q}}^{\frac{a+(k+1)b}{q}}f(t)dt-
  \frac{b}{q}\;f(\frac{a+kb}{q})\;\Big|\\
  \leq\;&\frac{b}{q}\;
  \|f\|_\infty +\sum_{k=0}^{q'-1}\;\frac{b}{q}\,(M_k-m_k) \leq
  (\|f\|_\infty+\Var(f))\,\frac{b}{q}\;.
\end{align*}
When $q'\neq 0$, since $|\,\omega_q(f)\,|\leq \|\omega_q\|\;
\|f\|_\infty$, we hence have
\begin{align*}
\overline{\omega_q}(f) &=  \int_0^1f(t)\,dt -\int_0^1f(t)\,dt+
\frac{b}{q}\,\omega_q(f) +\bigO\big(\frac{1}{q^2}\big)\,\omega_q(f)
\\ & =
\int_0^1f(t)\,dt +\bigO\Big(\frac{\|f\|_\infty+\Var(f)}{q}\Big)\;.
\end{align*}

For every $N\in\NN-\{0\}$,  with $N\geq a+b$, let us define
$$
\mu^-_N =\sum_{(m,\,n)\in I^-_N}\;\Delta_{\frac{m}{n}}\;\;
=\sum_{1\leq q\leq N,\; q\equiv a\!\!\!\mod b} \;\;\omega_q\;,
$$
which is a finitely supported measure on $[0,1]$. Its total mass is
equal to
$$
\|\mu^-_N\|=\sum_{1\leq q\leq N,\; q\equiv a\!\!\!\mod b} \|\omega_q\|=
\sum_{0\leq q'\leq \lfloor\frac{N-a}{b}\rfloor} (q'+\bigO(1))
=\frac{N^2}{2\,b^2}+\bigO(N)\;.
$$
Hence $\frac{1}{\|\mu^-_N\|} =\frac{2\,b^2}{N^2} + \bigO(\frac{1}{N^3})$.
For $f\in\BV([0,1])$, we have (taking $\|\omega_q\|\;
\overline{\omega_q}(f)=0$ if $q=a$)
\begin{align*}
  \frac{\mu^-_N(f)}{\|\mu^-_N\|}& =
  \frac{1}{\|\mu^-_N\|}\sum_{1\leq q\leq N,\; q\equiv a\!\!\!\mod b}
\;\;\|\omega_q\|\;\overline{\omega_q}(f)\\ &=
\int_0^1f(t)\,dt +\frac{1}{\|\mu^-_N\|}
\sum_{1\leq q\leq N,\; q\equiv a\!\!\!\mod b}
\bigO\big(\|f\|_\infty+\Var(f)\big)
\\ & =
\int_0^1f(t)\,dt +\bigO\Big(\frac{\|f\|_\infty+\Var(f)}{N}\Big)\;.
\end{align*}
Notice that $\ln$ is an increasing homeomorphism from $[0,1]$ to
$[-\infty,0]$.
%%%, and let us denote by $\ln_*$ the pushforwards of measures by $\ln$.
%For every bounded measurable function
%$f:[-\infty,0]\ra\RR$ with bounded variation, we hence have
%$\|f\circ\ln\|_\infty = \|f\|_\infty $ and $\Var(f\circ\ln)=\Var(f)$.
For every element $N\in\NN-\{0\}$, let us define
$$
\nu^\pm_N =\sum_{(m,\,n)\in I^\pm_N}\;\Delta_{\ln\frac{m}{n}}\;,
$$
so that $\nu^-_N=\ln_*\mu_N^-=\nu_N\mid_{]-\infty,0]}$, and
$\|\nu_N^-\|=\|\mu_N^-\|$. We have, for every $f\in BV(]-\infty,0])$,
\begin{align*}
  \frac{\nu^-_N(f)}{\|\nu^-_N\|}& =
  \frac{\mu^-_N(f\circ\ln)}{\|\mu^-_N\|} =\int_0^1f\circ\ln(t)\,dt
  +\bigO\Big(\frac{\|f\circ\ln\|_\infty+\Var(f\circ\ln)}{N}\Big)
  \\ & =\int_{-\infty}^0f(s)\;e^s\,ds +
  \bigO\Big(\frac{\|f\|_\infty+\Var(f)}{N}\Big)\;.
\end{align*}
Since $\nu_N=\nu^-_N+\nu^+_N$, since $\nu_N^+=\sg_*\nu^-_N$ 
%where $\sg:\RR\ra\RR$ is the change of sign map $t\mapsto -t$
, we have
$\|\nu^\pm_N\|=\frac{1}{2}\;\|\nu_N\|$ and the result follows.
\cqfd

\bigskip
Let us give some numerical illustrations of Theorem
\ref{theo:logpaircorrel} when $a=b=1$. For every
$N\in\NN-\{0\}$, let
$$
\D_N:s\mapsto
\frac{\card\{(p,q)\in I_N:\ln\frac pq\le s\}}{\card \;I_N}\;,
$$
which is the cumulative distribution function at time $N$ of the
differences of pairs of logarithms of integers, that is, for all
$s,s'\in \RR$ with $s<s'$, we have
$$
\frac{\nu_N}{\|\nu_N\|}\,(\,]s,s'\,]) =\D_N(s')-\D_N(s)\;.
$$
The above theorem says that the function $\D_N$ 
converges pointwise as $N\ra+\infty$ to the $C^1$ (but not $C^2$) function
$$
\D:s\mapsto\begin{cases}
\frac{1}{2}\;e^{s}&{\rm ~if~} s\leq 0\\
1-\frac{1}{2}\;e^{-s} &{\rm ~if~} s\geq 0\end {cases}
$$
(with derivative $\D'=g_{\L_\NN,1}$), which is the asymptotic cumulative
distribution function of the differences of pairs of logarithms of
integers. This is illustrated by the figure below, which shows
$\D_{15}$ in green.% and the limiting distribution $\D$ in blue.
%%%%%
%\todo{Is it the appropriate figure? If it should go to Section
%  \ref{sect:weightlogpaircorrel}, is there one possible here? Would
%  it be better looking by scaling it down vertically a bit?}
%%%%%

\begin{center}
\includegraphics[width=11cm]{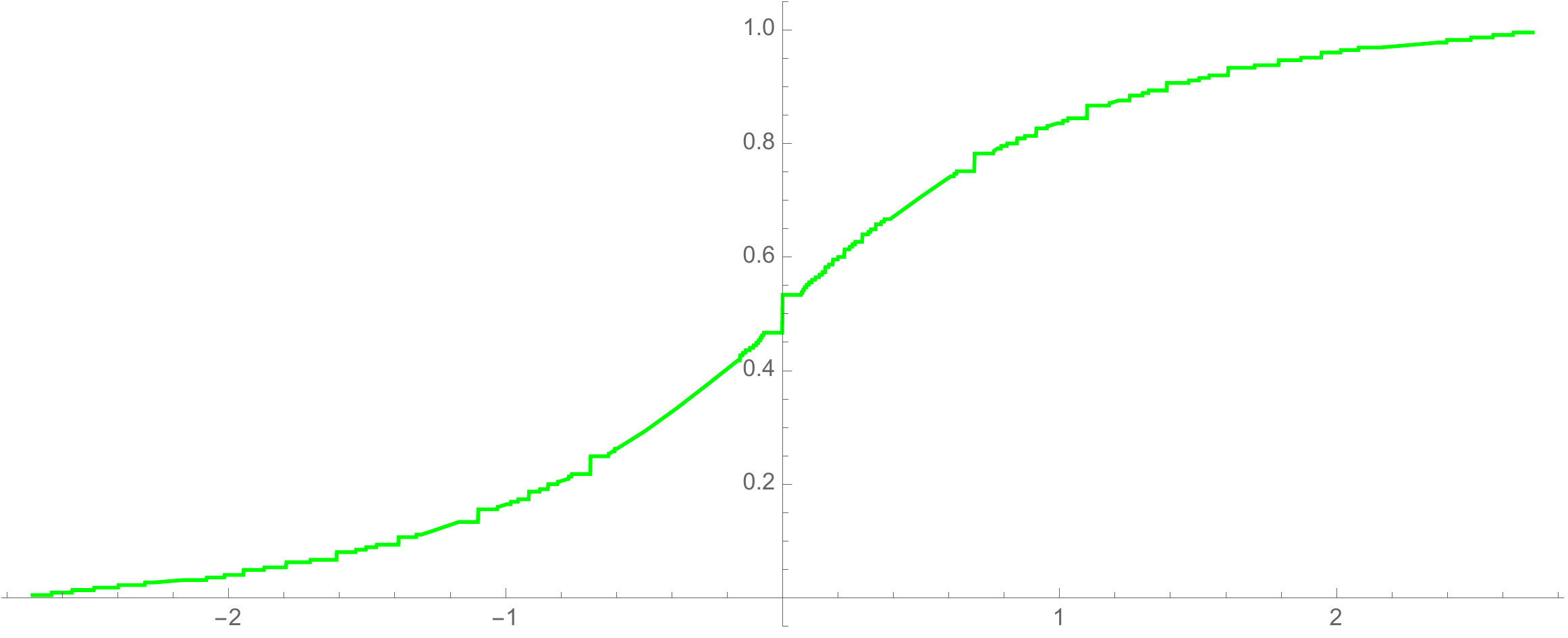}
\end{center}

\section{Pair correlations without weights and with scaling}
\label{sect:logpaircorrelscaled}

In this section, we study the pair correlations of logarithms of
integers at various scaling, now assumed to converge to $+\infty$. We
fix two nondecreasing positive functions, respectively $\psi:
\NN-\{0\}\ra \;]0,+\infty[$ and $\psi':\NN-\{0\}\ra \;]0,+\infty[$,
which will give the scaling factors on the difference of pairs of
logarithms and the renormalizing factors on their distribution.

For every $N\in\NN-\{0\}$, the (not normalised) {\em pair correlation
  measure} of the logarithms of integers congruent to $a$ modulo $b$
at time $N$ {\em with trivial multiplicities and with scaling}
$\psi(N)$ is the (Borel, positive) measure with finite support in
$\RR$ defined by
$$
\R^{\L_\NN,\psi}_N=\sum_{(m,\,n)\in I_N}\;\;\Delta_{\psi(N)(\ln m-\ln n)}\;,
$$
and the normalized one is $\frac{1}{\psi'(N)}\;\R^{\L_\NN,\psi}_N$.

\btheo\label{theo:logpaircorr} Assume that the nondecreasing positive
function $\psi$ satisfies ${\displaystyle \lim_{+\infty}}\;\psi
=+\infty$ and ${\displaystyle \lim_{N\ra+\infty}}\;\frac{\psi(N)}{N}
=\lambda_\psi\in [0,+\infty]$.  As $N\ra+\infty$, the measures
$\R^{\L_\NN,\psi}_N$ on $\RR$, normalized by $\psi'(N)$ as given
below, weak-star converge to a measure $g_{\L_\NN,\psi}\;\Leb_\RR$
absolutely continuous with respect to the Lebesgue measure on $\RR$,
$$
\frac{1}{\psi'(N)}\;\R^{\L_\NN,\psi}_N\;\;\;\weakstar\;\;\;
g_{\L_\NN,\psi}\;\Leb_\RR\;,
$$
with Radon-Nikodym derivative the function
$$ g_{\L_\NN,\psi}:t\mapsto\begin{cases} 0 &{\rm ~~if~~}
\lambda_\psi=+\infty, {\rm ~for~any~} \psi'\;, \\ \frac{1}{2\,b^2}
  &{\rm ~~if~~}\lambda_\psi=0{\rm ~and~}
  \psi'(N)=\frac{N^2}{\psi(N)}\;,\\ \frac{1}{2\,t^2}\;\big\lfloor
  \frac{|t|}{b\,\lambda_\psi}\big\rfloor \Big(\big\lfloor
  \frac{|t|}{b\,\lambda_\psi}\big\rfloor+1\Big) &{\rm
    ~~if~~}\lambda_\psi\neq 0,+\infty{\rm ~and~} \psi'(N)=\psi(N)\;.
\end{cases}
$$

Furthermore, if $\lambda_\psi\neq 0,+\infty$, for every $f\in C^1_{\rm
  c} (I)$ with support contained in $[-A,A]$, we have
$$
\frac{1}{\psi'(N)}\;\R^{\L_\NN,\psi}_N(f)=
\int_{s\in\RR}f(s)\;g_{\L_\NN,\psi}(s)\;ds+
\;\bigO_b\Big(A\,\|f\|_\infty
\big(\big|\lambda_\psi-\frac{\psi(N)}{N}\big|+\frac{A}{N}\big)+
\|f'\|_\infty\frac{A^3}{N}\Big)\,.
$$
\etheo

The pair correlation function $g_{\L_\NN,\psi}$ depends on $b$ but it
is independent of $a$.  The above result shows in particular that
renormalizing to probability measures (taking $\psi'(N)=N^2-N$) is
inappropriate, as the limiting measure would always be $0$.

When $\alpha=0$, $a=b=1$ and $\psi=\id^\alpha:N\ra N^\alpha$, the
measure $\frac{1}{\psi'(N)}\; \R^{\L_\NN,\psi}_N$ corresponds to the
one denoted by $\frac{1}{N^{2-\alpha}}\; \R^{\L_\NN,\,\id^\alpha}_N$
in the introduction. The above result thus implies the cases
$\alpha>0$ of Theorem \ref{theo:intro1} in the introduction, as well
as the comment about the scaling by the inverse of the average gap
$\psi(N)=\frac{N}{\ln N}$, for which $\lambda_\psi=0$.

The fact that $g_{\L_\NN,\psi}$ vanishes when $\lambda_\psi=+\infty$
means that the sequence of measures $\big(\frac{1}{\psi'(N)}\;
\R^{\L_\NN,\psi}_N \big)_{N\in\NN}$ on $\RR$ has a total loss of mass
at infinity. For error terms when $\lambda_\psi=+\infty$ and
$\lambda_\psi= 0$, see respectively Equation
\eqref{eq:errortermlambdaphiinfty} and Equation
\eqref{eq:errortermlambdaphizero}.

\medskip
\dem Note that the change of variables $(m,n)\mapsto (n,m)$ in $I_N$
proves that we have $\R^{\L_\NN,\psi}_N\!\mid_{]-\infty,0]}=\sg_*
\big(\R^{\L_\NN,\psi}_N\!  \mid_{[0,+\infty[}\big)$. We will thus
only study the convergence of the measures $\frac{1}{\psi'(N)}
\;\R^{\L_\NN,\psi}_N$ on $[0,+\infty[$, and deduce the global
result by the symmetry of $g_{\L_\NN,\psi}$.

\medskip
For every $N\in\NN-\{0\}$ and for every $p\in\NN$ with $p\equiv
0\!\!\mod b$ and $0<p<N$, let $N_p=\lfloor \frac{N-p-a}{b}\rfloor$,
let
\begin{equation}\label{eq:defiJsubpN}
J_{p,\,N}=\{q\in\NN :  1\leq q\leq N-p,\;q\equiv a\!\!\mod b\}=
\{a+kb : 0\leq k\leq N_p\}\;,
\end{equation}
and let
$$
\omega_{p,\,N}=\sum_{q\in J_{p,\,N}}\Delta_{\psi(N)\frac{p}{q}}
      {\rm ~~~and~~~}
\mu_N^+=\sum_{0<p<N,\;p\equiv 0\!\!\!\!\mod b}\omega_{p,\,N}\;.
$$
Then $\omega_{p,\,N}$ is a measure on $[0,+\infty[$, with finite
support contained in $[\frac{\psi(N)}{N-p}p,\psi(N)p]$, and with total
mass $\|\omega_{p,\,N}\| =N_p+1$. When this total mass is not $0$,
which occurs if and only if $N\geq p+a$, we denote by
$\overline{\omega_{p,\,N}} =\frac{\omega_{p,\,N}}
{\|\omega_{p,\,N}\|}$ the renormalisation of $\omega_{p,\,N}$ to a
probability measure on $]0,+\infty[\,$. The support of the measure
$\mu_N^+$ on $]0,+\infty[$ is contained in $[\frac{\psi(N)}{N},
\psi(N)N]$.  The motivation for the definition of the measure
$\mu_N^+$ comes from the following lemma.

\blemm\label{lem:relatRetmu}
For every $A>0$ and for every $f\in C^1_{\rm c}(\RR)$ with compact
support contained in $[0,A]$, we have, as $N\ra+\infty$,
$$
\big|\;\R^{\L_\NN,\psi}_N(f)-\mu_N^+(f)\,\big|=
\bigO\Big(A^3\;\|f'\|_\infty\;\big(\frac{N}{\psi(N)}\big)^2\Big)\;.
$$
\elemm

In particular, if $\frac{1}{\psi'(N)}\,\big(\frac{N}{\psi(N)}\big)^2$
tends to $0$ as $N\ra+\infty$, the measures $\frac{1}{\psi'(N)}
\,\R^{\L_\NN,\psi}_N\mid_{[0,+\infty]}$ and $\frac{1}{\psi'(N)}\,
\mu_N^+$ on $[0,+\infty]$ are asymptotic for the weak-star convergence
of measures on $\interval0{+\infty}$, and we will study the weak-star
convergence of the latter one.

\medskip
\dem By the change of variable $(p,q)\mapsto (m=p+q,n=q)$, we have
$$
\R^{\L_\NN,\psi}_N\mid_{[0,+\infty]}\;
=\sum_{(m,\,n)\in I_N^+}\;\;\Delta_{\psi(N)\ln\frac{m}{n}}
=\sum_{\substack{0<q\leq N-p,\;q\equiv a\!\!\!\mod b\\
    0<p<N,\;p\equiv 0\!\!\!\mod b}}
\;\;\Delta_{\psi(N)\ln(1+\frac{p}{q})}\;.
$$
By definition, we have
$$
\mu^+_N\;
=\sum_{\substack{0<q\leq N-p,\;q\equiv a\!\!\!\mod b\\
    0<p<N,\;p\equiv 0\!\!\!\mod b}}
\;\;\Delta_{\psi(N)\,\frac{p}{q}}\;.
$$
Since the support of $f$ is contained in $[0,A]$, if a pair $(p,q)$
occurs in the index of the sum defining either $\R^{\L_\NN,\psi}_N(f)$
or $\mu^+_N(f)$ with nonzero summand, then $\psi(N)\ln(1+\frac{p}{q})
\leq A$. This implies that $\frac{p}{q}=\bigO\big( \frac{A}{\psi(N)}
\big)$ since ${\displaystyle \lim_{+\infty}} \;\psi =+\infty$, and
that $p=\bigO\big(\frac{A\,N} {\psi(N)}\big)$ since $q\leq N$.  For
all $x,y\in [0,+\infty[\,$, we have
$$
|\Delta_x(f)-\Delta_y(f)|= |f(x)-f(y)|\leq \|f'\|_\infty|x-y|\;.
$$
Recall that $|\ln (1+t)-t|= \bigO(t^2)$ as $t\ra 0$. Hence, by a
uniform majoration of the terms of the sum below,
$$
\big|\;\R^{\L_\NN,\psi}_N(f)-\mu_N^+(f)\,\big|\leq
\sum_{\substack{1\leq q\leq N\\ 1\leq p\leq \bigO(\frac{A\,N}
    {\psi(N)})}} \|f'\|_\infty\;\psi(N)\;\bigO\Big((\frac{p}{q})^2\Big)
=\bigO\Big(A^3\;\|f'\|_\infty\;\big(\frac{N}{\psi(N)}\big)^2\Big)\;.
\;\;\;\Box
$$

\medskip
Let us now study the convergence properties of the (renormalized)
measures $\omega_{p,\,N}$ and of their sums $\mu^+_N$ as $N\ra+\infty$.

Let $\iota:\;]0,+\infty[\; \ra \;]0,+\infty[$ be the involutive
diffeomorphism $t\mapsto \frac{1}{t}$. We have
$$
\iota_*\omega_{p,\,N}=\sum_{q\,\in J_{p,\,N}}\Delta_{\frac{q}{\psi(N)p}}\;.
$$
As $q$ varies in $J_{p,\,N}$, the above Dirac masses are taken on the
distribution of points given by the following picture.

\begin{center}
\begin{picture}(0,0)%
\includegraphics{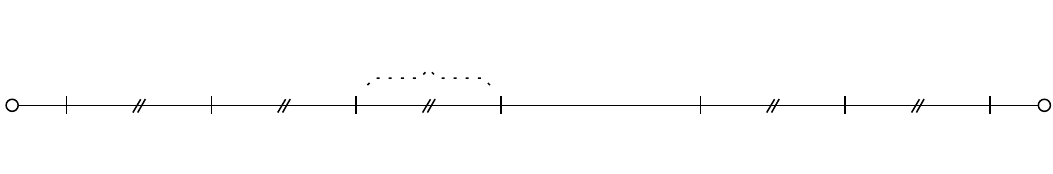}%
\end{picture}%
\setlength{\unitlength}{3812sp}%
\begingroup\makeatletter\ifx\SetFigFont\undefined%
\gdef\SetFigFont#1#2#3#4#5{%
  \reset@font\fontsize{#1}{#2pt}%
  \fontfamily{#3}\fontseries{#4}\fontshape{#5}%
  \selectfont}%
\fi\endgroup%
\begin{picture}(5228,858)(-59,-1390)
\put(1891,-691){\makebox(0,0)[lb]{\smash{{\SetFigFont{11}{13.2}{\rmdefault}{\mddefault}{\updefault}{\color[rgb]{0,0,0}$\frac{b}{\psi(N)p}$}%
}}}}
\put(5131,-871){\makebox(0,0)[lb]{\smash{{\SetFigFont{11}{13.2}{\rmdefault}{\mddefault}{\updefault}{\color[rgb]{0,0,0}$\frac{N-p}{\psi(N)p}$}%
}}}}
\put(3601,-1321){\makebox(0,0)[lb]{\smash{{\SetFigFont{11}{13.2}{\rmdefault}{\mddefault}{\updefault}{\color[rgb]{0,0,0}$\frac{a+(N_p-1)b}{\psi(N)p}$}%
}}}}
\put(4591,-1321){\makebox(0,0)[lb]{\smash{{\SetFigFont{11}{13.2}{\rmdefault}{\mddefault}{\updefault}{\color[rgb]{0,0,0}$\frac{a+N_pb}{\psi(N)p}$}%
}}}}
\put(  1,-1321){\makebox(0,0)[lb]{\smash{{\SetFigFont{11}{13.2}{\rmdefault}{\mddefault}{\updefault}{\color[rgb]{0,0,0}$\frac{a}{\psi(N)p}$}%
}}}}
\put(-44,-961){\makebox(0,0)[lb]{\smash{{\SetFigFont{11}{13.2}{\rmdefault}{\mddefault}{\updefault}{\color[rgb]{0,0,0}$0$}%
}}}}
\put(766,-1321){\makebox(0,0)[lb]{\smash{{\SetFigFont{11}{13.2}{\rmdefault}{\mddefault}{\updefault}{\color[rgb]{0,0,0}$\frac{a+b}{\psi(N)p}$}%
}}}}
\put(1531,-1321){\makebox(0,0)[lb]{\smash{{\SetFigFont{11}{13.2}{\rmdefault}{\mddefault}{\updefault}{\color[rgb]{0,0,0}$\frac{a+2b}{\psi(N)p}$}%
}}}}
\end{picture}%

\end{center}

As in the proof of Theorem \ref{theo:logpaircorrel}, for every $C^1$
function $f:\;]0,+\infty[\;\ra \RR$ with compact support, we have
\begin{align*}
  &\Big|\int_0^{\frac{N-p}{\psi(N)p}}f(t)\,dt-
  \frac{b}{\psi(N)p}\,\iota_*\omega_{p,\;N}(f)\,\Big|=
  \Big|\int_0^{\frac{N-p}{\psi(N)p}}f(t)\,dt- \sum_{q\,\in
    J_{p,\;N}}\;\frac{b}{\psi(N)p}\;f\big(\frac{q}{\psi(N)p}\big)
  \;\Big|\\ \leq\;&
  \Big|\int_0^{\frac{a}{\psi(N)p}}f(t)\,dt\Big|%\\ &\;\;\;
  +\sum_{k=0}^{N_p-1}\Big|
  \int_{\frac{a+kb}{\psi(N)p}}^{\frac{a+(k+1)b}{\psi(N)p}}f(t)dt-
   \;\frac{b}{\psi(N)p}\;f\Big(\frac{a+kb}{\psi(N)p}\Big)\;\Big| +
  \Big|\int_{\frac{a+N_pb}{\psi(N)p}}^{\frac{N-p}{\psi(N)p}}f(t)\,dt\,\Big|\\
  \leq\; &\big(2\,\|f\!\mid_{]0,\frac{N-p}{\psi(N)p}]}\|_\infty+
  \Var(f\!\mid_{]0,\frac{N-p}{\psi(N)p}]})\big)\frac{b}{\psi(N)p}\;.
\end{align*}
Hence for every $C^1$ function $f:\;]0,+\infty[\;\ra \RR$ with
    compact support, since $\iota$ is a diffeomorphism, we have
\begin{align*}
\omega_{p,\;N}(f)&=
\frac{\psi(N)p}{b}\;\int_0^{\frac{N-p}{\psi(N)p}} f\circ\iota(s)\;ds +
\bigO\big(\big\|f\circ\iota\!\mid_{]0,\frac{N-p}{\psi(N)p}]}\big\|_\infty+
\Var(f\circ\iota\!\mid_{]0,\frac{N-p}{\psi(N)p}]})\big)\\ & =
\frac{\psi(N)p}{b}\;
\int_{\frac{\psi(N)p}{N-p}}^{+\infty} f(t)\;\frac{dt}{t^2} +
\bigO\big(\big\|f\!\mid_{[\frac{\psi(N)p}{N-p},+\infty[}\big\|_\infty+
\Var(f\!\mid_{[\frac{\psi(N)p}{N-p},+\infty[})\big)\;.
\end{align*}
For every $t>0$, 
let $$\theta_N(t)=\frac{1}{t^2}
\sum_{\substack{0<p<N\\p\equiv 0\!\!\!\mod b}}\;\;\frac{p}{b}\;
\mathbbm{1}_{[\frac{\psi(N)p}{N-p},+\infty[}(t)\,.$$ Then 
\begin{align}
\theta_N(t)&= \frac{1}{t^2}\sum_{0<k<N/b}\;\;k\;
\mathbbm{1}_{[\frac{\psi(N)bk}{N-bk},+\infty[}(t)
=\frac{1}{t^2}\sum_{0<k\leq\frac{tN}{b(\psi(N)+t)}}\;\;k\nonumber\\ & =
\frac{1}{2\,t^2}\Big\lfloor \frac{tN}{b(\psi(N)+t)}\Big\rfloor
\Big(\Big\lfloor \frac{tN}{b(\psi(N)+t)}\Big\rfloor+1 \Big)\;.
\label{eq:valeurtheta}
\end{align}
In particular, we have $\theta_N(t)=0$ if and only if $t\in
[0,\frac{b\,\psi(N)}{N-b}[\,$.

Thus, if the support of $f$ is contained in the interval $[0,A]$,
since $\frac{\psi(N)p}{N-p}\leq A$ if and only if $p\leq
\frac{AN}{\psi(N)+A}$, we have,
\begin{align}
  \mu_N^+(f)&=\sum_{0<p<N,\;p\equiv 0\!\!\!\mod b}\omega_{p,\,N}(f)
  \nonumber\\ & = \psi(N)\int_0^{+\infty}f(t)
\Big(\sum_{\substack{0<p<N\\p\,\equiv\, 0\!\!\!\mod b}}
\frac{p}{b}\;\mathbbm{1}_{[\frac{\psi(N)p}{N-p},+\infty[}(t)\Big)
\;\frac{dt}{t^2} +\bigO\Big(\big(\|f\|_\infty+
\Var(f)\big)\frac{AN}{\psi(N)}\Big)  \nonumber\\ &
=\psi(N)\int_0^{+\infty}f(t)\;\theta_N(t)\;dt +
\bigO\Big(\big(\|f\|_\infty+\Var(f)\big)\frac{AN}{\psi(N)}\Big)
    \;.\label{eq:valeurfff}
\end{align}

\noindent{\bf Case 1. }
Assume first that $\lambda_\psi =+\infty$, that is, ${\displaystyle
  \lim_{N\ra+\infty}} \frac{N}{\psi(N)}=0$. Then for every $A\geq 1$,
if $N$ is large enough, then for every $t\in [0,A]$, we have
$\theta_N(t)=0$. Thus, whatever the normalizing function $\psi'$ is,
we have a total loss of mass at infinity~:
$$
\frac{1}{\psi'(N)}\;\mu^+_N\;\;\weakstar\;\; 0\;.
$$
More precisely, for every $C^1$ function $f:\;]0,+\infty[\;\ra
\RR$ with compact support contained in $[0,A]$, we have
$$
\frac{1}{\psi'(N)}\;\mu^+_N(f)=\bigO\Big(\big(\|f\|_\infty+
\Var(f)\big)\frac{AN}{\psi(N)\psi'(N)}\Big)\;.
$$
Since $\Var(f)=\int_0^A |f'| \leq A\;\|f'\|_\infty$, we have, by
Lemma \ref{lem:relatRetmu}, that $\frac{1}{\psi'(N)}\;\R^{\L_\NN,\psi}_N
\;\weakstar\; 0$ and
\begin{equation}\label{eq:errortermlambdaphiinfty}
\frac{1}{\psi'(N)}\;\R^{\L_\NN,\psi}_N(f)=\bigO\Big(\big(\|f\|_\infty+
\|f'\|_\infty\big)\frac{A^3N}{\psi(N)\psi'(N)}\Big)\;.
\end{equation}

\noindent{\bf Case 2. } Now assume that $\lambda_\psi= 0$, that is,
${\displaystyle \lim_{N\ra+\infty}} \frac{\psi(N)}{N}=0$. By Equation
\eqref{eq:valeurtheta}, if we have $t\geq \frac{b\psi(N)}{N-b}$, then
$$
\theta_N(t)=\frac{1}{2\,t^2}\big(\frac{tN}{b(\psi(N)+t)}+\bigO(1)\big)^2
= \frac{N^2}{2\,b^2\psi(N)^2}\Big(1 +\bigO\big(\frac{t}{\psi(N)}\big)+
\bigO\big(\frac{\psi(N)}{tN}\big) \Big)^2\;,
$$
therefore
\begin{equation}\label{eq:controlONt}
  \frac{\psi(N)^2}{N^2}\theta_N(t)=
   \frac{1}{2\,b^2} +\bigO\big(\frac{t}{\psi(N)}\big)+
\bigO\big(\frac{\psi(N)}{tN}\big) \;.
\end{equation}
Since $\theta_N$ vanishes on $[0,\frac{b\,\psi(N)}{N-b}\big[$, this
proves that $\frac{\psi(N)^2}{N^2}\;\theta_N$ is bounded on any
compact subset of $[0,+\infty[$, and pointwise converges to the
constant function $\frac{1}{2\,b^2}$. Hence by the Lebesgue dominated 
convergence theorem, we have
$$
\frac{\psi(N)}{N^2}\;\mu^+_N\;\;\weakstar\;\;
\frac{1}{2\,b^2}\;\Leb_{[0,+\infty[}\;.
$$

More precisely, for every $A\geq 3$, for every $C^1$ function $f:\;
]0,+\infty[\;\ra \RR$ with compact support contained in $[0,A]$, by
Equations \eqref{eq:valeurfff} and \eqref{eq:controlONt} and since
$\psi(N)\leq N$ for $N$ large enough, we have
\begin{align*}
  \frac{\psi(N)}{N^2}\;\mu_N^+(f) & =
  \frac{1}{2\,b^2}\int_{\frac{b\,\psi(N)}{N-b}}^{+\infty}f(t)\;dt
+\bigO\Big(\frac{1}{\psi(N)}\int_{\frac{b\,\psi(N)}{N-b}}^At\,|f(t)|\;dt\Big)
\\&\;\;\;\;+\bigO\Big(\frac{\psi(N)}{N}
\int_{\frac{b\,\psi(N)}{N-b}}^A\frac{1}{t}\,|f(t)|\;dt\Big)
+\bigO\Big(\big(\|f\|_\infty+\Var(f)\big)\frac{A}{N}\Big)
\\&=\frac{1}{2\,b^2}\int_{0}^{+\infty}f(t)\;dt
+\bigO\Big(\frac{\psi(N)}{N}\|f\|_\infty\Big)+
\bigO\Big(\frac{A^2}{\psi(N)}\|f\|_\infty\Big)
\\&\;\;\;\;+\bigO\Big(\frac{\psi(N)}{N}\,\|f\|_\infty
\big(\ln A-\ln\frac{b\,\psi(N)}{N-b}\big)\Big)
+\bigO\Big(\big(\|f\|_\infty+\|f'\|_\infty\big)\frac{A^2}{N}\Big)
\\&=\frac{1}{2\,b^2}\int_{0}^{+\infty}f(t)\;dt
+\bigO\Big(\|f\|_\infty\big(\frac{\psi(N)\ln A}{N}+
\frac{A^2}{\psi(N)}-\frac{\psi(N)}{N}\ln\frac{\psi(N)}{N}\big)\,\Big)
\\&\;\;\;\;+\bigO\Big(\|f'\|_\infty\frac{A^2}{N}\Big)
\;.
\end{align*}

Let $\psi'(N)=\frac{N^2}{\psi(N)}$. By
Lemma \ref{lem:relatRetmu}, we hence have
$$
\frac{1}{\psi'(N)}\;\R^{\L_\NN,\psi}_N
\;\;\weakstar\;\; \frac{1}{2\,b^2}\;\Leb_{\RR}
$$
Futhermore, for every $C^1$ function $f:\RR\ra \RR$ with
compact support contained in $[-A,A]$, we have
\begin{align}
  \frac{1}{\psi'(N)}\;\R^{\L_\NN,\psi}_N(f)&=
  \frac{1}{2\,b^2}\int_{\RR} f(t)\;dt
+\bigO\Big(\|f\|_\infty\big(\frac{\psi(N)\ln A}{N}+
\frac{A^2}{\psi(N)}-\frac{\psi(N)}{N}\ln\frac{\psi(N)}{N}\big)\,\Big)
\nonumber\\&\;\;\;\;\;\;\;\;\;\;\;\;\;\;\;\;\;\;\;\;\;\;\;\;\;\;\;\;
+\bigO\Big(\|f'\|_\infty\frac{A^3}{\psi(N)}\Big)\;.
\label{eq:errortermlambdaphizero}
\end{align}

\noindent{\bf Case 3. }
Let us finally assume that ${\displaystyle \lim_{N\ra+\infty}}
\frac{\psi(N)}{N}=\lambda_\psi$ belongs to $]0,+\infty[\,$.
Let us consider the map $\theta_\infty:[0,+\infty[\ra\RR$ defined by 
$$
t\mapsto\frac{1}{t^2}\sum_{k=1}^\infty\;\;k\;
\mathbbm{1}_{[b\lambda_\psi k,+\infty[}(t)=
\frac{1}{2\,t^2}\;\Big\lfloor \frac{t}{b\,\lambda_\psi}\Big\rfloor
\big(\Big\lfloor \frac{t}{b\,\lambda_\psi}\Big\rfloor+1\big)\;.
$$
It vanishes on $[0,b\,\lambda_\psi[\,$, is uniformly bounded, tends to
$\frac{1}{2\,b^2\,\lambda_\psi^2}$ as $t\ra+\infty$, and is
piecewise continous, with discontinuities at $b\,\lambda_\psi\NN-\{0\}$.
See the first picture in the introduction when $a=b=\lambda_\psi=1$.

By Equation \eqref{eq:valeurtheta}, the sequence of uniformly bounded
maps $(\theta_N)_{N\in\NN}$ converges almost everywhere to
$\theta_\infty$ (more precisely, it converges at least at every point
of $[0,+\infty[\;-\,b\,\lambda_\psi\NN\,$). Hence by Equation
    \eqref{eq:valeurfff} and by the Lebesgue dominated convergence
    theorem, we have
$$
\frac{1}{\psi(N)}\;\mu^+_N\;\;\weakstar\;\;
\theta_\infty\;\Leb_{[0,+\infty[}\;.
$$

Let $A\geq 1$ and $k\in\NN$. Note that $b\lambda_\psi k\leq A$ implies
that $k\leq \frac{A}{b\lambda_\psi}\leq \frac{2A}{b\lambda_\psi}$.  If
$N$ is large enough so that $\frac{\psi(N)}{N}\geq
\frac{\lambda_\psi}{2}$, then $\frac{\psi(N)bk}{N-bk}\leq A$ implies
that $k\leq \frac{AN}{b(\psi(N)+A)}\leq \frac{2A}{b\lambda_\psi}$.
Hence for every $t\in[0,A]$, we have
$$
|\;\theta_\infty(t)-\theta_N(t)\,|\leq \frac{1}{t^2}
\sum_{k=1}^{\frac{2A}{b\lambda_\psi}} k\;
\big|\,\mathbbm{1}_{[b\,\lambda_\psi k,+\infty[}(t)-
\mathbbm{1}_{[\frac{\psi(N)bk}{N-bk},+\infty[}(t)\,\big|\;.
$$
For every continuous function $f:[0,+\infty[\;\ra \RR$ with
compact support in $[0,A]$, we therefore have
\begin{align*}
\Big|\int_0^{+\infty} f\;(\theta_\infty-\theta_N)\;dt\,\Big|&=
\bigO\Big(\|f\|_\infty \sum_{k=1}^{\frac{2A}{b\lambda_\psi}} k\;
\big|\;b\,\lambda_\psi k-\frac{\psi(N)bk}{N-bk}\;\big|\; \Big)
\\&=\bigO\Big(A^3\;\|f\|_\infty \;\big|\,\lambda_\psi -
\frac{\psi(N)}{N}+\bigO(\frac{A}{N})\,\big|\, \Big)\;.
\end{align*}
By Equation \eqref{eq:valeurfff}, for every $C^1$ function
$f:[0,+\infty[\;\ra \RR$ with compact support in $[0,A]$, we thus have
\begin{align*}
\frac{1}{\psi(N)} \mu_N^+(f) &=\int f\;\theta_\infty\;d\Leb_{[0,+\infty[}
      \\ & \;\;\;\;+\bigO\Big(A^3\;\|f\|_\infty \;\big(\big|\,\lambda_\psi
      -\frac{\psi(N)}{N}\big|+\bigO(\frac{A}{N})\,\big)\; \Big)+
\bigO\Big(\|f'\|_\infty\frac{A^2}{N}\Big)\;.
\end{align*}

With $g_{\L_\NN,\psi}:\RR\ra\RR$ given by $t\mapsto \theta_\infty(|t|)$, by
Lemma \ref{lem:relatRetmu}, it follows that
$$
\frac{1}{\psi(N)}\;\R^{\L_\NN,\psi}_N
\;\;\weakstar\;\; g_{\L_\NN,\psi}\;\Leb_{\RR}\;.
$$
Futhermore, for every $A\geq 1$ and every $C^1$ function $f:\RR\ra
\RR$ with compact support contained in $[-A,A]$, we have
\begin{align}
  \frac{1}{\psi(N)}\;\R^{\L_\NN,\psi}_N(f)&=
  \int_{\RR} f(t)\;g_{\L_\NN,\psi}(t)\;dt
\\ & \;\;\;\;+\bigO\Big(A^3\;\|f\|_\infty\big( \;\big|\,\lambda_\psi
      -\frac{\psi(N)}{N}\big|+\frac{A}{N}\,\big)\; \Big)+
\bigO\Big(\|f'\|_\infty\frac{A^3}{N}\Big)\;.
\end{align}
This concludes the proof of Theorem \ref{theo:logpaircorr}.
\cqfd

\medskip
Let us give a numerical illustration of Theorem \ref{theo:logpaircorr}
when $a=b=1$ and $\psi(N)=N$.  The following figure shows in red an
approximation of the pair correlation function $g_{\L_\NN,\psi}$
computed using $\R^{\L_\NN,\psi}_{2000}$, and in blue the pair
correlation function $g_{\L_\NN,\psi}$ in the interval $[-4,4]$.

\begin{center}
\includegraphics[width=12cm]{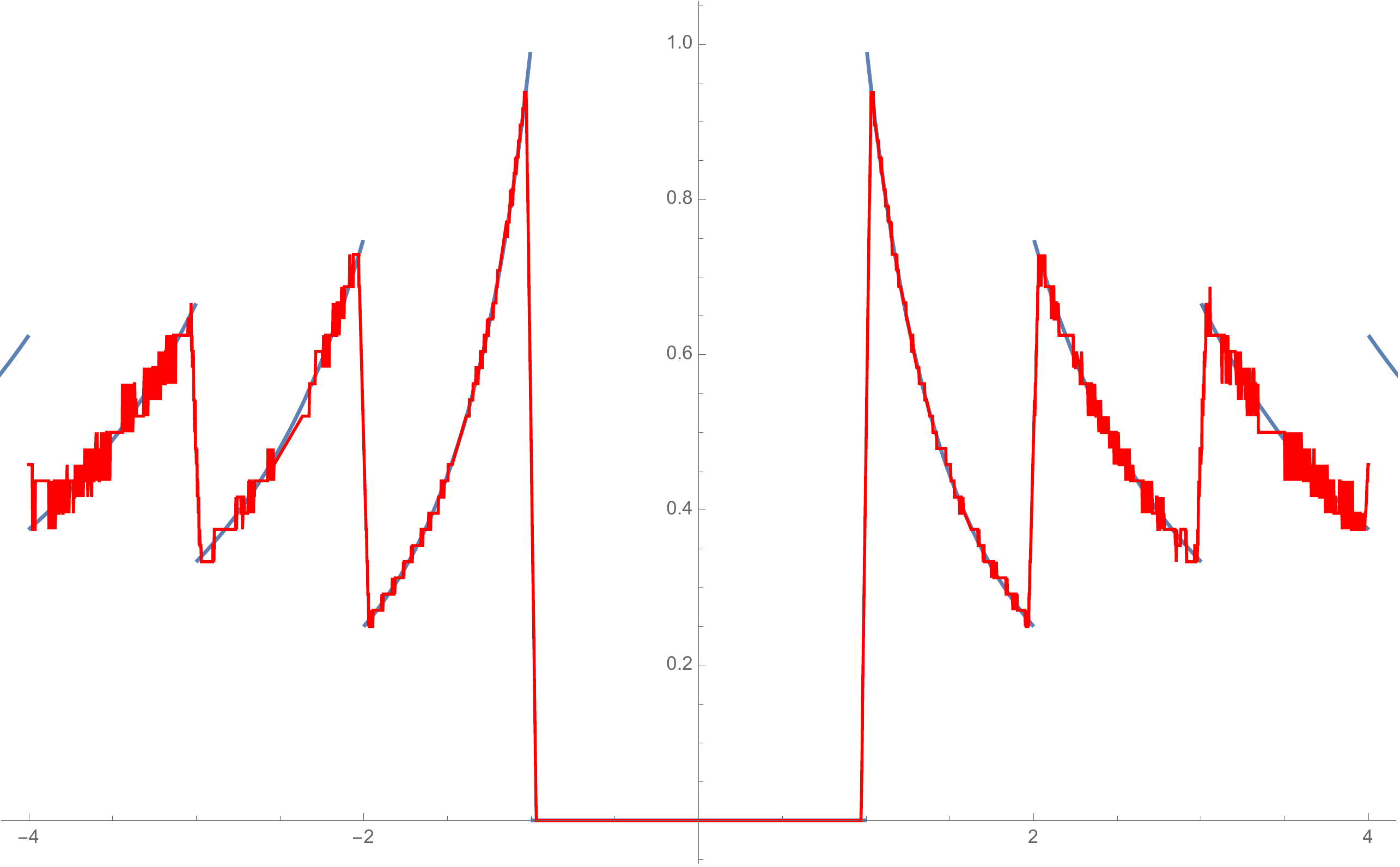}
\end{center}
%%
%\todo{change scale as for the figure above Proposition 5.5 ? Draw
%  horizontal line 1/2 ?}
%%

\section{Pair correlations with Euler weights without scaling}
\label{sect:weightlogpaircorrel}

In this section, we study the weighted family
$$
\L_\NN^{a,b,\varphi}=
\big(\,(L_N^{a,b})_{N\in\NN},\;\omega=\varphi\circ\exp\big)\,.
$$
The (not normalised) {\em pair correlation measure} of the logarithms
of integers congruent to $a$ modulo $b$ at time $N$ with
multiplicities given by the Euler function $\varphi$, for the trivial
scaling function, is
$$
\wt\nu_N =\sum_{(m,\,n)\in I_N}\;
\varphi(n)\;\varphi(m)\;\Delta_{\ln m-\ln n}\;.
$$
With the notation of the introduction, we have $\L_\NN^{1,1,\varphi}
=\L^\varphi_\NN$ and $\wt\nu_N= \R^{\L_\NN^{a,b,\varphi},1}_N$.

\btheo\label{theo:logpaircorrelphi} As $N\ra+\infty$, the measures
$\wt\nu_N$ on $\RR$, renormalized to be probability measures, weak-star
converge to the measure absolutely continuous with respect to the
Lebesgue measure on $\RR$, with Radon-Nikodym derivative the function
$g_{\L_\NN^{\varphi,a,b},1}:s\mapsto e^{-\,2\,|s|}$: 
$$
\frac{\wt\nu_N}{\|\,\wt\nu_N\|}\;\;\;\weakstar\;\;\;
g_{\L_\NN^{\varphi,a,b},1}\;\Leb_\RR\;.
$$
Furthermore, for all $f\in C^1_{\rm c}(\RR)$ and $\alpha\in
[\frac{1}{2},1[$, we have
$$
\frac{\wt\nu_N}{\|\,\wt\nu_N\|}(f)=\int_{s\in\RR}\;f(s)\;e^{-\,2\,|s|}\;ds
\;+\;\bigO_b\Big(\frac{\ln N}{N^{1-\alpha}}\,\|f\|_\infty+
\frac{1}{N^\alpha}\|e^{|s|}f'(s)\|_\infty\Big)\;.
$$
\etheo

When $a=b=1$, the measure $\wt\nu_N$ corresponds to the one denoted by
$\R^{\L_\NN^\varphi,1}_N$ in the introduction. The above result gives
the first assertion of Theorem \ref{theo:intro2} in the introduction,
with pair correlation function $g_{\L_\NN^\varphi,1}=
g_{\L_\NN^{\varphi,1,1},1}$.

\medskip
\dem The first assertion of Theorem \ref{theo:logpaircorrelphi}
follows from the second one, by taking for instance
$\alpha=\frac{1}{2}$ and by the density of $C^1$-smooth functions with
compact support in the space of continuous functions with compact support
on $\RR$.

For every $q\in\NN-\{0\}$ with $q\equiv a\!\!\mod b$, let
$q'\in\NN$ be such that $q=a+q'b$ and let $J_q$ be given by Equation
\eqref{eq:defJsubq}. We now define
$$
\wt\omega_q=\sum_{p\in J_q}\varphi(p)\;\Delta_{\frac{p}{q}}\;,
$$
which is a finitely supported measure on $[0,1]$, and nonzero if and
only if $q'\neq 0$. In order to compute its total mass, we will use
the following elementary adaptation of Mertens' formula (see for
example \cite[Thm. 330]{HarWri08}). We have not found its proof in the
literature, hence we provide one, due to Fouvry.

Let $(a,b)\in\NN-\{0\}$ be the greatest common divisor of $a$ and
$b$.  Let
$$
c_{a,b}=\frac{\varphi((a,b))}{b\,(a,b)}\;
\prod_{p {\rm ~prime},\; p\,\mid\,  b}(1-\frac{1}{p^2})^{-1}\;,
$$
and note that $c_{a,b}>0$ is uniformly bounded from above when $a$ and
$b$ vary in $\NN-\{0\}$, and tends to $0$ as $b\ra+\infty$. When
$a=b=1$, we have $c_{a,b}=1$, and the following result is exactly
Mertens' formula.

\blemm\label{lem:Fouvry1} There exists $C>0$ such that for all
integers $a,b\geq 1$ and real numbers $x\geq 1$, we have
$$
 \Big|\sum_{1\leq n\leq x,\;n\equiv a\!\!\!\mod b}\varphi(n)-
    \frac{3\,c_{a,b}}{\pi^2}\,x^2\;\Big|\leq C\;x\ln (2x)\;.
$$
\elemm

\dem (Fouvry) In this proof, for every function $g$ of a variable in
$[1,+\infty[\,$, possibly depending on parameters, we use the notation
$\bigO(g)$ in order to denote any function $f$ on $[1,+\infty[\,$
such that there exists a constant $C$, independent of the variable and
of all the parameters, such that $|f|\leq C|g|$. We do not need a more
precise error term.

Let $S(x,a,b)$ be the above sum. We refer for instance to
\cite{HarWri08} for the definition of the Möbius function
$\mu:\NN-\{0\}\ra \{-1,0,1\}$, of the Dirichlet convolution $f*g$ of
two maps $f,g:\NN-\{0\}\ra\RR$ and for the Möbius inversion formula,
which in particular gives that $\varphi=\mu*\id$. Hence
$$
S(x,a,b)
=\sum_{\substack{1\leq n\leq x\\n\equiv a\!\!\!\mod b}}
\;\sum_{md=n}\mu(d)\,m \;=\sum_{1\leq d\leq x}\mu(d)
\sum_{\substack{1\leq m\leq x/d\\md\equiv a\!\!\!\mod b}}m\;.
$$ Let us fix $d\geq 1$. Let us consider the congruence equation
$m\,d\equiv a\!\!\mod b$ with unknown $m$. It has no solution if the
greatest common divisor $(b,d)$ of $b$ and $d$ does not divide $a$. If
$(b,d)$ does divide $a$, let $a'=\frac{a}{(b,d)}$,
$b'=\frac{b}{(b,d)}$ and $d'=\frac{d}{(b,d)}$, so that the congruence
equation is equivalent to $m\,d'\equiv a'\!\!\mod b'$. Since $d'$ is
coprime with $b'$, it is invertible modulo $b'$, and we denote its
inverse by $\overline{d'}$. The congruence equation becomes $m\equiv
a'\;\overline{d'}\!\!\mod b'$.  The classical formula ${\displaystyle
\sum_{1\leq m\leq y, \;m\,\equiv\,a'\;\overline{d'}\!\!\!\mod b'}}1\;=
\frac{y}{b'}+\bigO(1)$ gives, by a summation by parts, the equality
$$
\sum_{1\leq m\leq y,\;m\,\equiv\, a'\;\overline{d'}\!\!\!\mod b'}m\;=
\frac{y^2}{2b'}+\bigO(y)\;.
$$
Therefore
\begin{align*}
S(x,a,b)&=\sum_{1\leq d\leq x,\;(b,d)\mid a}\mu(d)\Big(\frac{(b,d)}{2\,b}
\big(\frac{x}{d}\big)^2+\bigO\big(\frac{x}{d}\big)\Big)\\ &
=\frac{x^2}{2\,b}\Big(\sum_{1\leq d\leq x,\;(b,d)\mid a}
\mu(d)\frac{(b,d)}{d^2}\;\Big)+  \bigO(x\ln (2x))\;.
\end{align*}
Using the Eulerian product formula of the zeta function, giving
${\displaystyle\prod_{p {\rm ~prime}}}(1-\frac{1}{p^2})=
\frac{1}{\zeta(2)}= \frac{6}{\pi^2}$, and the expression
$\varphi(n)=n{\displaystyle\prod_{p{\rm~prime} ,\; p\,\mid\, n}}
(1-\frac{1}{p})$ of the Euler function in terms of the prime factors,
we have, by decomposing $d$ into prime powers and using the definition
of the Möbius function,
\begin{align*}
  \sum_{1\leq d\leq x,\;(b,d)\mid a}\mu(d)\frac{(b,d)}{d^2} &=
  \sum_{d\geq 1,\;(b,d)\mid a}\mu(d)\frac{(b,d)}{d^2} \;+
  \bigO\big(\frac{b}{x}\big)\\ &
  =\prod_{p {\rm ~prime},\; p\,\nmid\, b}(1-\frac{1}{p^2})
  \prod_{p {\rm ~prime},\; p\,\mid\,a,\;p\,\mid\, b}(1-\frac{1}{p})\;+
  \bigO\big(\frac{b}{x}\big)\\ &
  =\frac{1}{\zeta(2)}\prod_{p {\rm ~prime},\; p\,\mid\, b}
  (1-\frac{1}{p^2})^{-1}\;\frac{\varphi((a,b))}{(a,b)}+
  \bigO\big(\frac{b}{x}\big)=\frac{6\;b\;c_{a,b}}{\pi^2}+
  \bigO\big(\frac{b}{x}\big)\,.
\end{align*}
This proves Lemma \ref{lem:Fouvry1}.
\cqfd

\medskip
Lemma \ref{lem:Fouvry1} says that if $q'\neq 0$ (that is, when $q> a$), then
\begin{equation}\label{eq:totmasswtomsubq}
\|\,\wt\omega_q\|=\frac{3\,c_{a,b}\,q^2}{\pi^2}+\bigO(q\ln (2q))\;,
\end{equation}
and in particular $\frac{1}{\|\,\wt\omega_q\|}=
\frac{\pi^2}{3\,c_{a,b}\,q^2}\big(1+\bigO(\frac{\ln q}{q})\big)$.

\blemm \label{lem:calcwtomega} We have as $q\ra+\infty$,
$$
\frac{\wt\omega_q}{\|\,\wt\omega_q\|}\;\;\weakstar
\;\;2\,t\;d\Leb_{[0,1]}(t)\;.
$$
More precisely, for all $f\in C^1([0,1])$ and $\alpha\in
[\frac{1}{2},1\,[\,$, we have
$$
\frac{\wt\omega_q}{\|\,\wt\omega_q\|}(f)=
\int_0^1\;2\,t\,f(t)\;\;dt
\;+\;\bigO\Big(\frac{\ln q}{q^{1-\alpha}}\|f\|_\infty+
\frac{1}{q^\alpha}\|f'\|_\infty\Big)\;.
$$    
\elemm

\dem The first assertion follows from the second one, by taking for
instance $\alpha=\frac{1}{2}$ and by the density of $C^1$-smooth
functions in the space of continuous functions on $[0,1]$.

Let $Q=\lfloor q^\alpha\rfloor\in\NN-\{0\}$. For all $n\in\{0,\dots,
Q-1\}$ and $t\in\;]\,n\,q^{-\alpha},(n+1)q^{-\alpha}]$, we have by the
mean value theorem
$$
f(t)=f(n\,q^{-\alpha})+\bigO(q^{-\alpha}\,\|f'\|_\infty)\;.
$$
Since $n\leq Q\leq q^\alpha$, we hence have
\begin{align}
  \int_{n\,q^{-\alpha}}^{(n+1)q^{-\alpha}}2\,t\,f(t)\;dt&=
  \big(f(n\,q^{-\alpha})+\bigO(q^{-\alpha}\,\|f'\|_\infty)\big)
  \int_{n\,q^{-\alpha}}^{(n+1)q^{-\alpha}}2\,t\;dt \nonumber 
\\ & =\frac{1}{q^\alpha}\big((2\,n+1)\,q^{-\alpha} f(n\,q^{-\alpha}) +
\bigO(n\,q^{-2\alpha}\|f'\|_\infty)\,\big)\;.
\label{eq:controlint}
\end{align}
Using twice Equation \eqref{eq:totmasswtomsubq} and the formula for
$\frac{1}{\|\,\wt\omega_q\|}$ following it, we have
\begin{align}
  &\sum_{p\,\in\,]\,n\,q^{1-\alpha},\,(n+1)q^{1-\alpha}]\cap J_q}
f\big(\frac{p}{q}\big)\;\frac{1}{\|\,\wt\omega_q\|}\;\varphi(p)\;
\nonumber\\ =\; &
  \big(f(n\,q^{-\alpha})+\bigO(q^{-\alpha}\,\|f'\|_\infty )\big)
  \frac{\pi^2}{3\,c_{a,b}\,q^2}\,
  \big(1+\bigO\big(\frac{\ln q}{q}\big)\big)
  \nonumber\\  &\;\;\;\times
  \Big(\,\frac{3\,c_{a,b}((n+1)^2-n^2)}{\pi^2}q^{2-2\alpha}+
  \bigO\big((n+1)q^{1-\alpha}\ln((n+1)q^{1-\alpha})\big)\,\Big)
  \nonumber\\ =\; & \frac{1}{q^\alpha}
  \Big((2\,n+1)\,q^{-\alpha}\,f(n\,q^{-\alpha})+
  \bigO\big(\,\frac{n}{q^{2\alpha}} \,\|f'\|_\infty+
  \frac{n\ln q}{q}\,\|f\|_\infty\big)\Big)\;.
\label{eq:controlmoyphi}
\end{align}
Again using Equation \eqref{eq:totmasswtomsubq}, we have
\begin{equation}\label{eq:controlrestsum}
\sum_{p\,\in\,]\,Q\,q^{1-\alpha},\,q[\cap J_q}
\frac{\varphi(p)}{\|\,\wt\omega_q\|}\;\big|f\big(\frac{p}{q}\big)\big|=
\bigO\big(\frac{q^2-(Q\,q^{1-\alpha})^2}{q^2}\,\|f\|_\infty\big)
=\bigO\big(q^{-\alpha}\,\|f\|_\infty\big)\;.
\end{equation}
By cutting the sum defining $\wt\omega_q$ and the integral from $0$ to
$1$, by using \eqref{eq:controlint}, \eqref{eq:controlmoyphi} and
\eqref{eq:controlrestsum}, since $Q\leq q^\alpha$ and $|1-
Q\,q^{-\alpha}| =\bigO(q^{-\alpha})$, we have (using $\alpha\geq
\frac{1}{2}$ for the last equality)
\begin{align*}
  &\Big|\;\frac{\wt\omega_q}{\|\,\wt\omega_q\|}(f)-
  \int_0^1\;2\,t\,f(t)\;\;dt\;\Big|
\\=\; & \Big|\;\sum_{n=0}^{Q-1} \Big(
  \sum_{p\,\in\,]\,n\,q^{1-\alpha},\,(n+1)q^{1-\alpha}]\cap J_q}
\frac{\varphi(p)}{\|\,\wt\omega_q\|}\;f\big(\frac{p}{q}\big)-
\int_{n\,q^{-\alpha}}^{(n+1)q^{-\alpha}}2\,t\,f(t)\;dt\;\Big) \;\Big|
\\ &\;\;+
\sum_{p\,\in\,]\,Q\,q^{1-\alpha},\,q[\cap J_q}
\frac{\varphi(p)}{\|\,\wt\omega_q\|}\;
\big|f\big(\frac{p}{q}\big)\big|\;\;+\;
\int_{Q\,q^{-\alpha}}^12\,t\,|f(t)|\;dt\;
\\=\; &\bigO(q^{-\alpha}\,\|f\|_\infty)+\frac{1}{q^\alpha}\sum_{n=0}^{Q-1}
\bigO\big(\,\frac{n}{q^{2\alpha}}\|f'\|_\infty
+\frac{n\ln q}{q}\|f\|_\infty+n\,q^{-2\alpha}\,\|f'\|_\infty)\big)
\\=\; &\bigO(q^{-\alpha}\,\|f\|_\infty)+\frac{1}{q^\alpha}
\bigO\big(\,\frac{Q^2}{q^{2\alpha}}\|f'\|_\infty
+\frac{Q^2\ln q}{q}\|f\|_\infty)\big)
\\=\; &
\bigO\big(\frac{\ln q}{q^{1-\alpha}}\,\|f\|_\infty+
\frac{\|f'\|_\infty}{q^\alpha}\big)\;.
\end{align*}
This proves Lemma \ref{lem:calcwtomega}.
\cqfd

\medskip
For every $N\in\NN-\{0\}$, let us define
$$
\wt\mu^-_N =\sum_{(m,\,n)\in I^-_N}\;\varphi(m)\;\varphi(n)\;
\Delta_{\frac{m}{n}}\;\;
=\sum_{1\leq q\leq N,\;q\equiv a\!\!\!\mod b}\;\varphi(q)\;\wt\omega_q\;,
$$
which is a finitely supported measure on $[0,1]$. By Lemma
\ref{lem:Fouvry1}, its total mass is
\begin{align*}
  \|\wt\mu^-_N\|&=
  \sum_{1\leq q\leq N,\;q\equiv a\!\!\!\mod b}\;\varphi(q)\;\|\wt\omega_q\|
  =\sum_{(m,n)\in I_N^-} \;\varphi(m)\;\varphi(n)\\ &=
\frac{1}{2}\Big(
\big(\sum_{\substack{1\leq q\leq N\\q\equiv a\!\!\!\mod b}} \varphi(q)\big)^2-
\sum_{\substack{1\leq q\leq N\\q\equiv a\!\!\!\mod b}}\varphi(q)^2\Big)
=\frac{9\,c_{a,b}^2\,N^4}{2\,\pi^4}+\bigO(N^3\ln N)\;.
\end{align*}
For $f\in C^1([0,1])$, by Lemma \ref{lem:calcwtomega}, since
$q^{1+\alpha}(\ln q)\leq N^{1+\alpha}(\ln N)$ and $q^{2-\alpha}\leq
N^{2-\alpha}$ when $q$ occurs in the summations below, we have
\begin{align*}
  \frac{\wt\mu^-_N(f)}{\|\wt\mu^-_N\|}& = \frac{1}{\|\wt\mu^-_N\|}
  \sum_{1\leq q\leq N,\;q\equiv a\!\!\!\mod b}\;\varphi(q)
\;\|\wt\omega_q\|\;\frac{\wt\omega_q(f)}{\|\wt\omega_q\|}\\ &=
\int_0^12\,t\,f(t)\,dt +\frac{1}{\|\wt\mu^-_N\|}
\sum_{\substack{1\leq q\leq N\\q\equiv a\!\!\!\mod b}}
\;\varphi(q)\;\|\wt\omega_q\|\,
\bigO\big(\frac{\ln q}{q^{1-\alpha}}\,\|f\|_\infty+
\frac{1}{q^\alpha}\|f'\|_\infty\big)
\\ & =
\int_0^12\,t\,f(t)\,dt +\bigO\big(\frac{1}{N^{4}}
\sum_{\substack{1\leq q\leq N\\q\equiv a\!\!\!\mod b}}\;\varphi(q)
\big(q^{1+\alpha}(\ln q)\,\|f\|_\infty+q^{2-\alpha}\,\|f'\|_\infty\big)
\\ & =\int_0^12\,t\,f(t)\,dt +\bigO\big(\frac{\ln N}{N^{1-\alpha}}
\,\|f\|_\infty+ \frac{1}{N^\alpha}\|f'\|_\infty\big)\;.
\end{align*}

For every $N\in\NN-\{0\}$, let us define
$$
\wt\nu^\pm_N =\sum_{(m,\,n)\in I^\pm_N}\;\varphi(m)\;\varphi(n)\;
\Delta_{\ln\frac{m}{n}}\;,
$$
so that $\wt\nu^-_N=\ln_*\wt\mu_N^-=\wt\nu_N\!\mid_{]-\infty,\,0]}$,
and $\|\wt\nu_N^-\|=\|\wt\mu_N^-\|$. We have, for every $f\in C^1_{\rm
  c} (]-\infty,0])$,
\begin{align*}
  \frac{\wt\nu^-_N(f)}{\|\wt\nu^-_N\|}& =
  \frac{\wt\mu^-_N(f\circ\ln)}{\|\wt\mu^-_N\|} =
  \int_0^12\,t\,f\circ\ln(t)\,dt
  +\bigO\big(\frac{\ln N}{N^{1-\alpha}}\,\|f\circ\ln\|_\infty+
\frac{1}{N^\alpha}\|(f\circ\ln)'\|_\infty\big)
\\ & =\int_{-\infty}^02\,f(s)\;e^{2s}\,ds +
\bigO\big(\frac{\ln N}{N^{1-\alpha}}\,\|f\|_\infty+
\frac{1}{N^\alpha}\|e^{-s}f'(s)\|_\infty\big)\;.
\end{align*}
Since $\wt\nu_N=\wt\nu^-_N+\wt\nu^+_N$, since $\wt\nu_N^+=
\sg_*\wt\nu^-_N$, we have $\|\wt\nu^\pm_N\|=\frac{1}{2}\;\|\wt\nu_N\|$
and Theorem \ref{theo:logpaircorrelphi} follows.
\cqfd

Let us give some numerical illustrations of Theorem
\ref{theo:logpaircorrelphi} with $a=b=1$. For every $N\in\NN-\{0\}$,
let
$$
\wt\D_N:s\mapsto
\frac{\sum_{0\leq m\neq n\leq N \;:\; (\ln m-\ln n)\le s} \varphi(m)\;\varphi(n)}
{\sum_{0\leq m\neq n\leq N}\varphi(m)\;\varphi(n)}\;.
$$
This is the cumulative distribution function at time $N$ of the
unscaled differences of the logarithms of natural numbers weighted by
the Euler function, that is, for all $s,s'\in \RR$ with $s<s'$, we
have
$$
\frac{\R^{\L_{\NN}^\varphi,\,1}_N}{\|\R^{\L_{\NN}^\varphi,\,1}_N\|}\,(\,]s,s'\,]) =
\wt\D_N(s')-\wt\D_N(s)\;.
$$
Theorem \ref{theo:logpaircorrelphi} with $a=b=1$ says that the
function $\wt\D_N$ pointwise converges as $N\ra+\infty$ to the $C^1$
(but not $C^2$) function
$$
\wt\D:s\mapsto\Big\{\begin{array}{l}
\frac{1}{2}\;e^{2s}{\rm ~if~} s\leq 0\\
1-\frac{1}{2}\;e^{-2s} {\rm ~if~} s\geq 0\end {array}
$$
(with derivative $\wt\D'=g_{\L_\NN^{\varphi},1}$).  This is
illustrated by the figure below, which shows $\wt\D_{15}$ in red.

\begin{center}
\includegraphics[width=11cm]{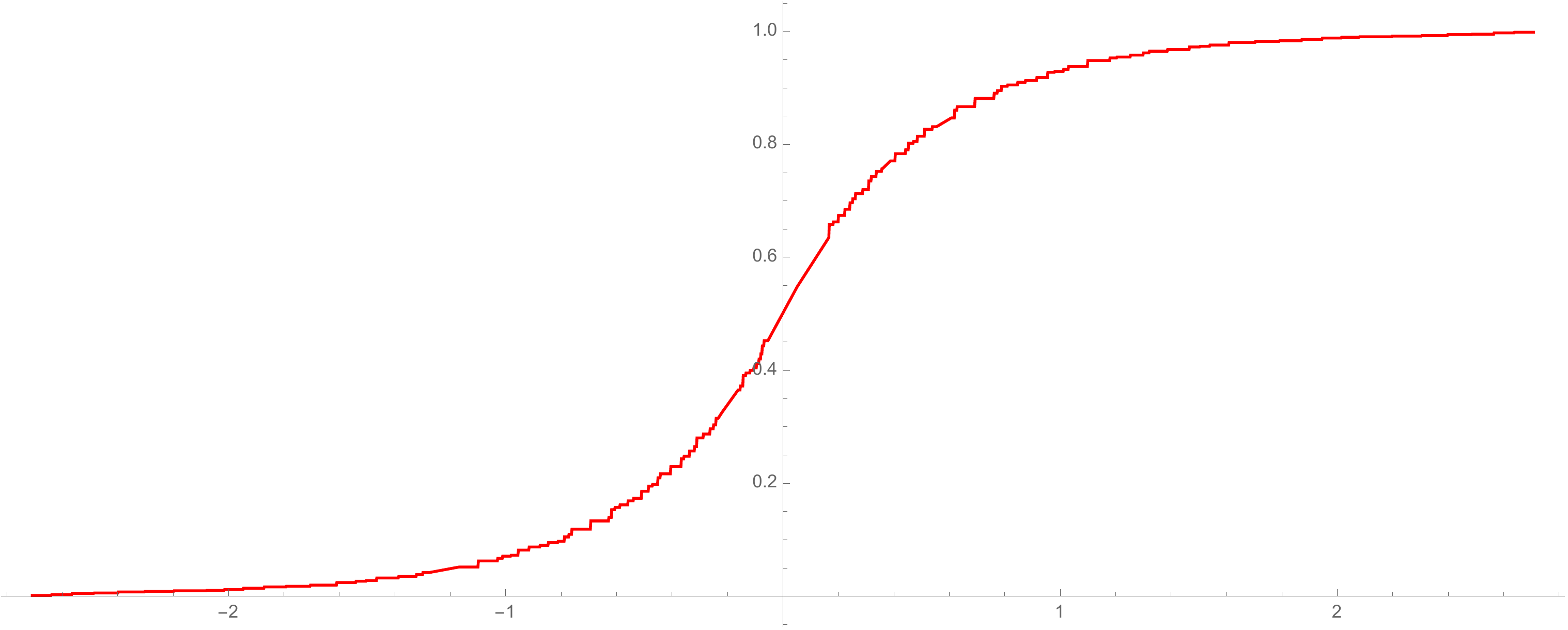}
\end{center}

\section{Pair correlations with Euler weights and linear scaling}
\label{sect:weightscaledlogpaircorrel}

In this section, we study the pair correlations of the family
$\L_\NN^{\varphi,a,b}$ defined at the beginning of Section
\ref{sect:weightlogpaircorrel}, now with a linear scaling. We leave to
the reader the study of a general scaling $\psi$, assumed to converge
to $+\infty$.  For every $N\in\NN-\{0\}$, the (not normalised) {\em
  pair correlation measure} of the logarithms of integers, congruent
to $a$ modulo $b$, at time $N$ with multiplicities given by the Euler
function and with scaling $N$ is the (Borel, positive) measure with
finite support in $\RR$ defined by
$$
\wt\R_N=\sum_{(m,\,n)\in I_N}\;\varphi(m)\;\varphi(n)\;\Delta_{N(\ln m-\ln n)}\;.
$$
With the notation of the introduction, we have
$\wt\R_N=\R_N^{\L_\NN^{\varphi,a,b},\;\id}$.

For every $k\in\NN-\{0\}$, we consider the arithmetic constant
$c_{a,b,k}$ defined in Equation \eqref{eq:cabk} of Appendix
\ref{appendixFouvry}. Note that $c_{a,b,k}>0$ is uniformly bounded
from above when $a,b,k$ vary in $\NN-\{0\}$, and has a positive lower
bound on $a$ and $k$ when $b$ is fixed, by Equation \eqref{eq:cabkc}
of Appendix \ref{appendixFouvry}.

\btheo\label{theo:logpaircorrelphipsi}
As $N\ra+\infty$, the family $\big(\frac{1}{N^3}\;
\wt\R_N\big)_{N\in\NN}$ of measures on $\RR$ weak-star converges to the
measure absolutely continuous with respect to the Lebesgue measure on
$\RR$, with Radon-Nikodym derivative the function
$$
g_{\L_\NN^{\varphi,a,b},\,\id}:s\mapsto \frac{1}{s^4}
\sum_{1\leq k\leq\lfloor\,|s|\,\rfloor,\;k\equiv 0\!\!\!\mod b}
c_{a,b,k} \;k^3\;,
$$
that is, as $N\ra+\infty$,
$$
\frac{1}{N^3}\;\wt\R_N\;\;\;\weakstar\;\;\;
g_{\L_\NN^{\varphi,a,b},\,\id}\;\Leb_\RR\;.
$$
Furthermore, for all $f\in C^1(\RR)$ with compact support contained in
$[-A,A]$, where $A\geq 1$, and for any $\alpha\in\; [\frac{1}{2},
  1\,[\,$, we have
$$
\frac{1}{N^3}\;\wt\R_N(f)=
\int_{s\in\RR}\;f(s)\;g_{\L_\NN^{\varphi,a,b},\,\id}(s)\;ds
\;+\;\bigO_b\Big(\frac{A^7\ln^2N}{N^{1-\alpha}}\,\|f\|_\infty+
\frac{A^3}{N^{\alpha}}\|f'\|_\infty\Big)\;.
$$
\etheo

When $a=b=1$, the measure $\wt\R_N$ corresponds to the one denoted by
$\R^{\L_\NN^\varphi,\,\id}_N$ in the introduction. The above result
gives the second assertion of Theorem \ref{theo:intro2} in the
introduction, with pair correlation function $g_{\L_\NN^\varphi,\,\id}
=g_{\L_\NN^{\varphi,1,1},\,\id}$, using Mirsky's value of $c_{1,1,k}$
given by Equation \eqref{eq:c11kmirsky}, as explained in Remark
\ref{rem:retrouvMirsky} of Appendix \ref{appendixFouvry}.

Note that, as the proof below shows, the total mass of
$\wt\R_N$ is equivalent to $c\,N^4$ as $N\ra+\infty$, for some
constant $c>0$, hence renormalising $\wt\R_N$ to be a probability
measure makes it weak-star converge to the zero measure on the
noncompact space $\RR$ (a total loss of mass phenomenon).

\medskip
\dem The first assertion of Theorem \ref{theo:logpaircorrelphipsi}
follows from the second one, by taking for instance $\alpha=
\frac{1}{2}$ and by the density of $C^1$-smooth functions with compact
support in the space of continuous functions with compact support on
$\RR$.

The change of variables $(m,n)\mapsto(n,m)$ in $I_N$ gives
$\wt\R_N|_{]-\infty,0]}=\sg_* \big(\wt\R_N|_{[0,+\infty[}\big)$. We
will thus only study the convergence of the measures $\frac{1}{N^3}\;
\wt\R_N$ on $[0,+\infty[$, and deduce the global result by the
symmetry of $g_{\L_\NN^{\varphi,a,b},\,\id}$.
    
\medskip
For every $N\in\NN-\{0\}$ and for every $p\in\NN$ with $p\equiv
0\!\!\mod b$ and $0<p<N$, let $J_{p,\,N}$ be given par Equation
\eqref{eq:defiJsubpN}. We now define the key auxiliary measure by
$$
\wt\omega_{p,\,N}=
\sum_{q\in J_{p,\,N}}\varphi(q)\;\varphi(q+p)\;\Delta_{\frac{q}{N\,p}}\;.
$$
Then $\wt\omega_{p,\,N}$ is a measure on $[0,+\infty[$, with finite
support contained in $[\frac{1}{N\,p},\frac{N-p}{N\,p}]$, hence in
$[0,1]$. The measure $\wt\omega_{p,\,N}$ is nonzero if and only if
$N\geq a+p$. In order to compute its total mass, we use an adaptation
with congruences of a formula by Mirsky (see \cite[Thm.~9]{Mirsky49})
proved in the appendix by Fouvry. Theorem \ref{theo:mirskyfouvry} says
that if $N\geq a+p$, then
\begin{equation}\label{eq:totmasswtomsubpN}
  \|\,\wt\omega_{p,N}\|=\frac{c_{a,b,p}}{3}\,(N-p)^3+
  \bigO(N^2(p+\ln(N-p))^2)\;,
\end{equation}
and in particular $\frac{1}{\|\,\wt\omega_{p,N}\|}=
\frac{3}{c_{a,b,p}\,(N-p)^3} \big(1+
\bigO(\frac{N^2(p+\ln(N-p))^2}{(N-p)^3})\big)$.

\medskip
The next result implies that the measures $\wt\omega_{p,N}$, once
normalized to be probability measures, weak-star converge to the
measure $d\mu(t)=3\big(\frac{N\,p}{N-p}\big)^3t^2\;
d\Leb_{[\frac{1}{N\,p},\frac{N-p}{N\,p}]}$, which is absolutely continuous
with respect to the Lebesgue measure on the interval $[\frac{1}{N\,p},
  \frac{N-p}{N\,p}]$.

\blemm\label{lem:distribomegasubpN} For every $p\in\NN$ with $0<p<N$
and $p\equiv 0\!\!\mod b$, for every $\alpha\in \;]\,0,1\,[\,$ and for
    every $f\in C^1_{\rm c}(]0,1])$, we have
\begin{align*}
\frac{\wt\omega_{p,N}}{\|\wt\omega_{p,N}\|}(f)&=
\int_{\frac{1}{N\,p}}^{\frac{N-p}{N\,p}}3
\big(\frac{N\,p}{N-p}\big)^3\;t^2\;f(t)\;dt\\ &\;\;\;\; +
\bigO\Big(\big(\frac{(N\,p)^3}{(N-p)^3(N\,p)^{\alpha}}+
\frac{(N\,p)^6(p+\ln(N\,p))^2}{(N-p)^6(N\,p)^{1-\alpha}}\big)\;\|f\|_\infty+
\frac{1}{(N\,p)^{\alpha}}\;\|f'\|_\infty\Big)\;.
\end{align*}
\elemm

\dem As in the proof of Lemma \ref{lem:calcwtomega}, we will
estimate the difference of the main terms in the above centered
formula by cutting the sum defining the renormalized measure
$\wt\omega_{p,N}$ and by cutting similarly the integral from
$\frac{1}{N\,p}$ to $\frac{N-p}{N\,p}$.

Let $Q=\lfloor (N\,p)^\alpha\frac{N-p}{N\,p}\rfloor\in\NN$. For all
$n\in\{0,\dots, Q-1\}$, we thus define
$$
S_n=
\sum_{q\,\in\,]\,n\,(N\,p)^{1-\alpha},\,(n+1)(N\,p)^{1-\alpha}]\cap J_{p,N}}
f\big(\frac{q}{N\,p}\big)\;\frac{1}{\|\,\wt\omega_{p,N}\|}
\;\varphi(q)\;\varphi(q+p)
$$
and
$$
I_n=\int_{n\,(N\,p)^{-\alpha}}^{(n+1)(N\,p)^{-\alpha}}
3\big(\frac{N\,p}{N-p}\big)^3\;t^2\;f(t)\;dt\;.
$$
Let us also define the following remaining terms
$$
S_{\rm end}=
\sum_{q\,\in\,]\,Q\,(N\,p)^{1-\alpha},\,N-p]\cap J_{p,N}}
f\big(\frac{q}{N\,p}\big)\;\frac{1}{\|\,\wt\omega_{p,N}\|}
\;\varphi(q)\;\varphi(q+p)
$$
and
$$
I_{\rm end}=\int_{Q\,(N\,p)^{-\alpha}}^{\frac{N-p}{N\,p}}
3\big(\frac{N\,p}{N-p}\big)^3\;t^2\;f(t)\;dt\;.
$$

For all $t\in\;]\,n\,(N\,p)^{-\alpha},(n+1)(N\,p)^{-\alpha}]$, we
have by the mean value theorem
$$
f(t)=f(n\,(N\,p)^{-\alpha})+\bigO((N\,p)^{-\alpha}\,\|f'\|_\infty)\;.
$$
Using twice Theorem \ref{theo:mirskyfouvry} and the formula for
$\frac{1}{\|\,\wt\omega_{p,N}\|}$ following Equation
\eqref{eq:totmasswtomsubpN}, and using the inequality $(n+1)\leq
Q\leq(N\,p)^{\alpha}\frac{N-p}{N\,p}$, we have
\begin{align}
S_n&=\Big(f(n\,(N\,p)^{-\alpha})+\bigO((N\,p)^{-\alpha}\,
\|f'\|_\infty)\Big)\Big(\frac{1}{(N-p)^3}
\big(1+\bigO\big(\frac{N^2(p+\ln(N-p))^2}{(N-p)^3}\big)\,\big)\Big)
\nonumber\\  &\;\;\;\;\;\times
\Big(\big((n+1)(N\,p)^{1-\alpha}\big)^3-\big(n(N\,p)^{1-\alpha}\big)^3
\nonumber\\  &\;\;\;\;\;\;\;\;\;\;\;+
\bigO\big(\,\big(p+(n+1)(N\,p)^{1-\alpha}\big)^2\big(p+
\ln((n+1)(N\,p)^{1-\alpha})\big)^2\,\big)\Big)
\nonumber\\  =& \frac{(N\,p)^3}{(N-p)^3(N\,p)^{\alpha}}
\Big((3n^2+3n+1)(N\,p)^{-2\alpha}f(n\,(N\,p)^{-\alpha}) +
\bigO\big(\frac{n^2}{(N\,p)^{3\alpha}}\;\|f'\|_\infty\big)
\nonumber\\  &\;\;\;\;\;\;\;\;\;\;\;
+\bigO\big(\,\big(\,\frac{n^2N^2(p+\ln(N-p))^2}{(N-p)^3(N\,p)^{2\alpha}}
+\frac{(n+\frac{p}{(Np)^{1-\alpha}})^2(p+\ln(N\,p))^2}
{(N\,p)^{1+\alpha}}\,\big)\,\|f\|_\infty\,\big)\Big)\;.
\label{eq:controlmoyphipsi}
\end{align}
We also have
\begin{align}
  I_n&=\Big(f(n\,(N\,p)^{-\alpha})+\bigO((N\,p)^{-\alpha}\,\|f'\|_\infty)\Big)
  \int_{n\,(N\,p)^{-\alpha}}^{(n+1)(N\,p)^{-\alpha}}
  3\big(\frac{N\,p}{N-p}\big)^3\;t^2\;dt\nonumber\\ & =
  \frac{(N\,p)^3}{(N-p)^3(N\,p)^{\alpha}}
  \Big((3n^2+3n+1)(N\,p)^{-2\alpha}f(n\,(N\,p)^{-\alpha}) +
  \bigO\big(\frac{n^2}{(N\,p)^{3\alpha}}\;\|f'\|_\infty\big)\Big)\;.
\label{eq:controlintpsi}
\end{align}
Similarly, since $Q\geq (N\,p)^\alpha\frac{N-p}{N\,p}-1$ and $N-p\leq
N\,p$, we have
\begin{equation}\label{eq:controlSend}
S_{\rm end}=\bigO\Big(\frac{(N-p)^3-(Q(N\,p)^{1-\alpha})^3}{(N-p)^3}
\;\|f\|_\infty\Big) =
\bigO\Big(\frac{(N\,p)^3}{(N-p)^3(N\,p)^{\alpha}}\;\|f\|_\infty\Big)
\end{equation}
and
\begin{equation}\label{eq:controlIend}
  I_{\rm end}=\bigO\Big(\big(\frac{N\,p}{N-p}\big)^3
  \Big(\big(\frac{N-p}{N\,p}\big)^3-(Q(N\,p)^{-\alpha})^3\Big)
\;\|f\|_\infty\Big) =
\bigO\Big(\frac{(N\,p)^3}{(N-p)^3(N\,p)^{\alpha}}\;\|f\|_\infty\Big)\;.
\end{equation}
Note that
${\displaystyle\sum_{n=0}^{Q-1}}\big(n+\frac{p}{(Np)^{1-\alpha}}\big)^2
=\bigO\big(\big(Q+\frac{p}{(Np)^{1-\alpha}}\big)^3\big)=
\bigO\big(\frac{(N\,p)^{3\alpha}}{p^3}\big)$. Putting
together Equations \eqref{eq:controlmoyphipsi},
\eqref{eq:controlintpsi}, \eqref{eq:controlSend} and
\eqref{eq:controlIend}, we have (again using $N-p \leq N\,p$)
\begin{align*}
 &\Big|\;\frac{\wt\omega_{p,N}}{\|\wt\omega_{p,N}\|}(f)-
  \int_{\frac{1}{N\,p}}^{\frac{N-p}{N\,p}}
  3\big(\frac{N\,p}{N-p}\big)^3\;t^2\;f(t)\;dt\;\Big|=
  \Big|\sum_{n=0}^{Q-1} (S_n-I_n) + S_{\rm end}-I_{\rm end}\;\Big|
  \\\leq \;&\sum_{n=0}^{Q-1} |S_n-I_n| + |S_{\rm end}|+|I_{\rm end}|
  \\=\;&
  \bigO\Big(\frac{(N\,p)^3}{(N-p)^3(N\,p)^{\alpha}}\;\|f\|_\infty\Big)+
  \frac{(N\,p)^3}{(N-p)^3(N\,p)^{\alpha}}\sum_{n=0}^{Q-1}
  \bigO\Big(\;\frac{n^2}{(N\,p)^{3\alpha}}\;\|f'\|_\infty+
  \\ &\;\;\;\;\;\;\;
  \big(\,\frac{n^2N^2(p+\ln(N-p))^2}{(N-p)^3(N\,p)^{2\alpha}}
  +\frac{\big(n+\frac{p}{(Np)^{1-\alpha}}\big)^2(p+\ln(N\,p))^2}
  {(N\,p)^{1+\alpha}}\,\big)\,
  \|f\|_\infty+\frac{n^2}{(N\,p)^{3\alpha}}\;\|f'\|_\infty\Big) \\=\;&
  \bigO\Big(\frac{(N\,p)^3}{(N-p)^3(N\,p)^{\alpha}}\;\|f\|_\infty+
  \frac{1}{(N\,p)^{\alpha}}\;\|f'\|_\infty+
  \frac{(N\,p)^6(p+\ln(N\,p))^2}{(N-p)^6(N\,p)^{1-\alpha}}
  \,\|f\|_\infty\Big)\;.
\end{align*}
This proves Lemma \ref{lem:distribomegasubpN}.
\cqfd

\medskip
Now, let us introduce the sum
$$
\wt\mu^+_N=\sum_{0<p<N,\;p\equiv0\!\!\!\mod b}\;\iota_*\wt\omega_{p,N}\;
=\sum_{\substack{1\leq q\leq N-p,\;0<p<N\\q\equiv a\!\!\!\mod b,\;p\equiv 0\!\!\!\mod b}}
\;\varphi(q)\;\varphi(q+p)\;
\Delta_{N\frac{p}{q}}\;,
$$
where as previously $\iota:t\mapsto \frac{1}{t}$ (noting that the
measures $\wt\omega_{p,\,N}$ are supported in $]0,+\infty[\,$).

\blemm\label{lem:comparRsubNmusubN} For every $f\in C^1([0,+\infty[)$
with compact support contained in $[0,A]$, where $A\geq 1$, we have
$$
\big|\;\wt \R_N(f)-\wt\mu_N(f)\,\big|=\bigO(A^3\,N^2\,\|f'\|_\infty)\;.
$$
\elemm

\dem Using the change of variables $(p,q)\mapsto (m=p+q,n=q)$, we have
$$
\wt\R_N\!\mid_{[0,+\infty[}=\sum_{(m,n)\in I_N^+}\;\varphi(m)\;\varphi(n)\;
\Delta_{N\ln\frac{m}{n}}\;
=\sum_{\substack{0<p<N,\;1\leq q\leq N-p\\
    p\equiv 0\!\!\!\mod b,\;q\equiv a\!\!\!\mod b}}
\;\varphi(q)\;\varphi(q+p)\;
\Delta_{N \ln(1+\frac{p}{q})}\;.
$$
As in the proof of Lemma \ref{lem:relatRetmu}, since the support of
$f$ is contained in $[0,A]$, if a pair $(p,q)$ occurs in the index of
the sum defining either $\wt\R_N(f)$ or $\wt\mu^+_N(f)$ with nonzero
corresponding summand, then $\frac{p}{q}=\bigO\big(\frac{A}{N}\big)$ and
$p=\bigO(A)$. By the mean value theorem, we then have
\begin{align*}
  \big|\,f\big(N\,\frac{p}{q}\big)-
  f\big(N\ln(1+\frac{p}{q})\big)\,\big|&\leq
  \|f'\|_\infty\big|\,N\,\frac{p}{q}-N\ln(1+\frac{p}{q})\,\big|\\ &
  =\|f'\|_\infty\,N\,\bigO\big(\big(\frac{p}{q}\big)^2\big)=
  \bigO\big(\frac{A^2}{N}\|f'\|_\infty\big)\;.
\end{align*}
Thus, using Theorem \ref{theo:mirskyfouvry} and Equation
\eqref{eq:cabkc} in Appendix \ref{appendixFouvry}, we have
\begin{align*}
  \big|\,\wt\R_N(f)-\wt\mu^+_N(f)\big|&\leq
  \sum_{1\leq p\leq \bigO(A),\;1\leq q\leq N}
  \;\varphi(q)\;\varphi(q+p)
  \bigO\big(\frac{A^2}{N}\|f'\|_\infty\big)\\ &\leq
  \sum_{\substack{1\leq p\leq \bigO(A)}} \bigO(N^3)
  \bigO\big(\frac{A^2}{N}\|f'\|_\infty\big)=
  \bigO\big(A^3\,N^2\,\|f'\|_\infty\big)\;.
\end{align*}
This proves Lemma \ref{lem:comparRsubNmusubN}.
\cqfd

\blemm\label{lem:distribmusubN} For all $\alpha\in[\frac{1}{2},1[$
and $f\in C^1([0,+\infty[)$ with compact support contained in $[0,A]$,
where $A\geq 1$, we have, as $N\ra+\infty$,
$$%\begin{align*}
\frac{1}{N^3}\;\wt\mu_N(f)=
\int_0^\infty f(s)\;g_{\L_\NN^{\varphi,a,b},\,\id}(s)\;ds+
\bigO\Big(\frac{A^7\ln^2N}{N^{1-\alpha}}\;\|f\|_\infty+
\frac{A^3}{N^\alpha}\;\|f'\|_\infty\Big)\;.
$$%\end{align*}
\elemm

\dem Let $A$ and $f$ be as in the statement.  Since the support of
$\wt\omega_{p,\,N}$ is contained in $[\frac{1}{N\,p}, \frac{N-p}
{N\,p}]$, the support of $ \iota_*\wt\omega_{p,\,N}$ is contained in
$[\frac{N\,p}{N-p},N\,p]$.  In particular the measures $\wt\mu_N$ and
$g_{\L_\NN^{\varphi,a,b},\,\id}(s) \;ds$ both vanish on $[0,1]$. Hence
we may assume that the support of $f$ is contained in $[1,+\infty[\,$,
    so that the support of $f\circ\iota$ is contained in $]0,1]$.

Note that $\|f\circ\iota\|_\infty= \|f\|_\infty$ and
$\|(f\circ\iota)'\|_\infty= \|t^2f'(t)\|_\infty\leq
A^2\|f'\|_\infty$, since the support of $f'$ is contained in $[0,A]$.

By the definition of $\wt\mu^+_N $, by Equation
\eqref{eq:totmasswtomsubpN} and Lemma \ref{lem:distribomegasubpN},
since $N-p\leq N$ and by the restriction on $p$, explained in the
proof of Lemma \ref{lem:comparRsubNmusubN}, in the summation defining
$\wt\mu^+_N(f) $ due to the support of $f$, we have
\begin{align*}
  \wt\mu^+_N(f)&=\sum_{0<p<N,\;p\equiv0\!\!\!\mod b}\;
  \iota_*\wt\omega_{p,\,N}(f)=
  \sum_{0<p<N,\;p\equiv0\!\!\!\mod b}\;
  \|\wt\omega_{p,\,N}\|\;\frac{\wt\omega_{p,N}}{\|\wt\omega_{p,\,N}\|}
  (f\circ \iota)\\ &
  =\sum_{0<p<N,\;p\equiv0\!\!\!\mod b,\; p\leq\bigO(A)}\;
  \Big(\frac{c_{a,b,p}}{3}\,(N-p)^3+
  \bigO\big(N^2(p+\ln(N-p))^2\big)\Big)
  \\ &\;\;\;\times \Big(\int_{\frac{1}{N\,p}}^{\frac{N-p}{N\,p}}3
\big(\frac{N\,p}{N-p}\big)^3\;t^2\;f\circ \iota(t)\;dt\\ &\;\;\; +
\bigO\Big(\big(\frac{(N\,p)^3}{(N-p)^3(N\,p)^{\alpha}}+
\frac{(N\,p)^6(p+\ln(N\,p))^2}{(N-p)^6(N\,p)^{1-\alpha}}\big)\;
\|f\circ \iota\|_\infty+
\frac{1}{(N\,p)^{\alpha}}\;\|(f\circ \iota)'\|_\infty\Big)\Big)
\\ &=N^3\sum_{0<p<N,\;p\equiv0\!\!\!\mod b,\; p\leq\bigO(A)}c_{a,b,p}\;p^3
\int_{\frac{N\,p}{N-p}}^{N\,p}\frac{1}{s^4}\;f(s)\;ds
\\ &\;\;\;\; +N^3\sum_{0<p<N,\;p\equiv0\!\!\!\mod b,\; p\leq\bigO(A)}
\bigO\big(\;\frac{p^3N^2(p+\ln(N-p))^2}{(N-p)^3}\;\|f\|_\infty\big)
\\ &\;\;\;\;\;\;\;\;\;\;\;\;\;\;\;\;\;\;\;\;\;\;\;\;\;\;\;
+\bigO\Big(\big(\frac{p^{3-\alpha}}{N^{\alpha}}+
\frac{N^3p^{5+\alpha}(p+\ln(N\,p))^2}{(N-p)^3N^{1-\alpha}}\big)\;
\|f\|_\infty+\frac{A^2}{(N\,p)^{\alpha}}\;\|f'\|_\infty\Big)\Big)\;.
\end{align*}
Noting that
${\displaystyle\int_p^{\frac{N\,p}{N-p}}}\frac{1}{s^4}\;f(s)=
  \bigO\big(\frac{N^2p+N\,p^2+p^3}{(N\,p)^3}\;\|f\|_\infty\big)$,
  we therefore have, for $N$ large compared to $A$,
\begin{align*}
  \frac{1}{N^3}\;\wt\mu^+_N(f)&=
  \Big(\sum_{0<p<N,\;p\equiv0\!\!\!\mod b}c_{a,b,p}\;p^3
  \int_{p}^{+\infty}\frac{1}{s^4}\;f(s)\;ds\Big)+
  \bigO\big(\frac{A^4}{N}\;\|f\|_\infty\big) \\ &\;\;\;\;\;+
  \bigO\Big(\frac{A^4\ln^2(N)}{N}\;\|f\|_\infty+
  \big(\frac{A^4}{N^{\alpha}}+\frac{A^7\ln^2(N)}{N^{1-\alpha}}\big)\;
  \|f\|_\infty+\frac{A^3}{N^{\alpha}}\;\|f'\|_\infty\Big)
  \\ &=\int_0^{+\infty}f(s)\;g_{\L_\NN^{\varphi,a,b},\,\id}(s)\;ds+
  \bigO\Big(\;\big(\frac{A^4}{N^{\alpha}}+
  \frac{A^7\ln^2(N)}{N^{1-\alpha}}\big)\;
  \|f\|_\infty+\frac{A^3}{N^{\alpha}}\;\|f'\|_\infty\Big)\;.
\end{align*}
This proves Lemma \ref{lem:distribmusubN}, using that
$\frac{A^4}{N^{\alpha}}=\bigO(\frac{A^7\ln^2(N)}{N^{1-\alpha}})$ if
$\alpha\geq \frac{1}{2}$.
\cqfd

\medskip
Theorem \ref{theo:logpaircorrelphipsi} now follows from Lemmas
\ref{lem:comparRsubNmusubN} and \ref{lem:distribmusubN}, as explained
in the beginning of the proof.  \cqfd

\medskip
We close this section with a numerical illustration of Theorem
\ref{theo:logpaircorrelphipsi} when $a=b=1$.  The following figure
shows in red an approximation of the pair correlation function
$g_{\L_\NN^\varphi,\,\id}=g_{\L_\NN^{\varphi,1,1},\,\id}$ computed
using $\wt\R_{2000}$ in the interval $[-10,10]$, to be compared with
the graph of $g_{\L_\NN^\varphi,\,\id}$ in the introduction.

\begin{center}
\includegraphics[width=12cm]{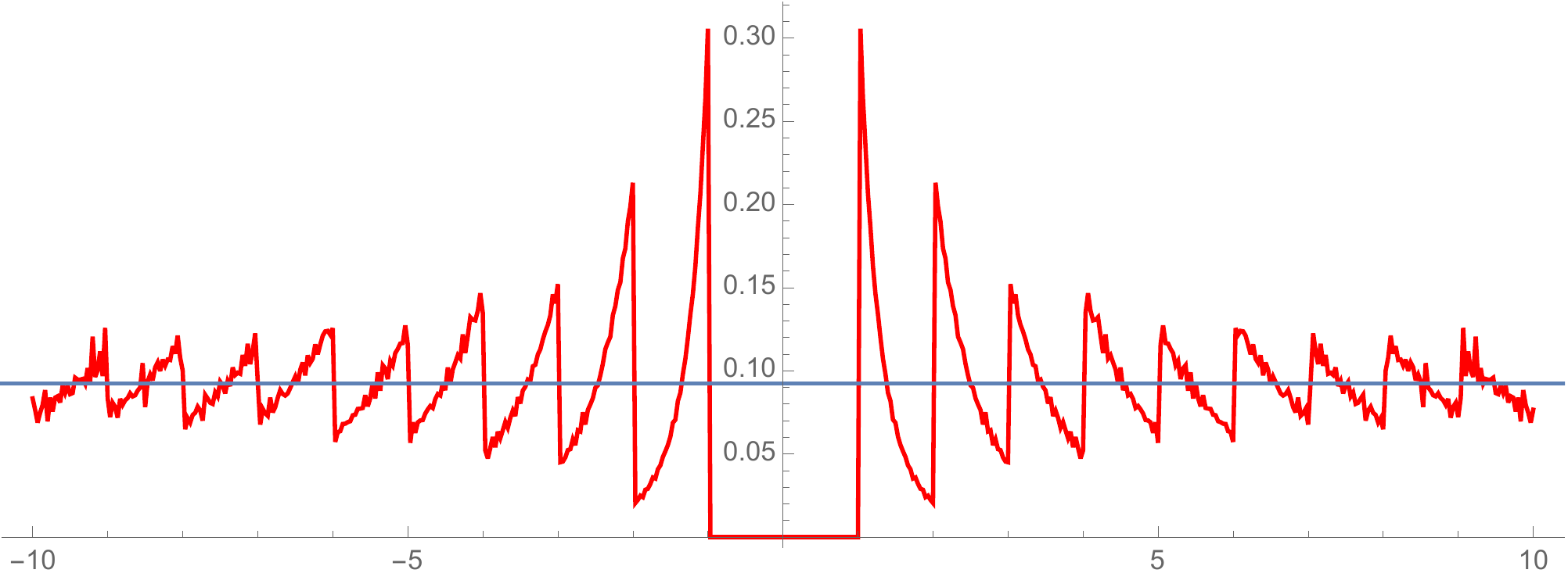}
\end{center}

The fact that the graph of $g_{\L_\NN^\varphi,\,\id}$ has a
horizontal asymptote near $\pm\infty$ follows from the following
result.

\bprop\label{prop:Fouvry2} We have $\lim_{s\ra\pm\infty}\;
g_{\L_\NN^\varphi,\,\id}(s) =
\frac{1}{4}\;{\displaystyle\prod_{p {\rm ~prime}}}\big(1-\frac{2}{p^2}\big)
\big(1+\frac{1}{p^2(p^2-2)}\big)$.
\eprop

\dem (Fouvry) In this proof, we use the same convention concerning
$\bigO(\cdot)$ as in the beginning of the proof of Lemma
\ref{lem:Fouvry1}.

We consider the multiplicative\footnote{Recall that an
  arithmetic function $f$ is {\it multiplicative} if $f(1)=1$ and for
  all coprime integers $m,n$, we have $f(mn)=f(m)f(n)$.} function
$f:n\mapsto \!\displaystyle{\prod_{p {\rm ~prime},\;p\,\mid\, n}}\!
\big(1+\frac{1}{p(p^2-2)}\big)$ and the constant $C_1=
\displaystyle{\prod_{p {\rm ~prime}}}\big(1+\frac{1} {p^2(p^2-2)}\big)$. Let us
prove that uniformly in $x\geq
1$, we have 
\begin{equation}\label{eq:asympcube}
\sum_{1\leq n\leq x}n^3f(n)= \frac{C_1}{4}\;x^4\;+\bigO(x^3)\;.
\end{equation}
By Equation \eqref{eq:gLvarphiid}, this proves Proposition
\ref{prop:Fouvry2}.

Let $g=f*\mu$ be the Dirichlet convolution of $f$ with the Möbius
function $\mu$. Then $g$ is multiplicative. For every prime $p$, we
have
$$
g(p)=f(p)\,\mu(1)+f(1)\,\mu(p)=\frac{1}{p(p^2-2)}
$$
and $g(p^k)=f(p^k)\,\mu(1)+f(p^{k-1})\,\mu(p)=0$ for every $k\geq 2$.
Therefore, for every $m\geq 1$, we have
$$
g(m)=\mu(m)^2\;\prod_{p {\rm ~prime},\; p\,\mid\, m}\;\frac{1}{p(p^2-2)}\;.
$$

\blemm\label{lem:controlpointg} For every $m\geq 1$, we have $0\leq
g(m)\leq m^{-3}\displaystyle{\prod_{p {\rm ~prime}}}
\big(1-\frac{2}{p^2}\big)^{-1}$.
\elemm

\dem 
This is immediate if $\mu(m)=0$. Otherwise, $m=p_1\dots p_k$ with
$p_1,\dots ,p_k$ pairwise distinct primes, and
$$
0\leq m^{3}g(m)=\prod_{i=1}^k \frac{p_i^3}{p_i(p_i^2-2)}=
\prod_{i=1}^k \big(1-\frac{2}{p_i^2}\big)^{-1}\leq
\prod_{p {\rm ~prime}}\big(1-\frac{2}{p^2}\big)^{-1}<+\infty\;. \;\;\;\Box
$$

Therefore, using the Möbius inversion formula $f=g*{\bf 1}$ for the
first equality, Lemma \ref{lem:controlpointg} for the fifth equality
and an Eulerian product (since $g$ is multiplicative and vanishes on
integers divisible by a nontrivial square) for the sixth equality, we
have, with $S(x)= \sum_{1\leq k\leq x} f(k)$,
\begin{align*} S(x)&
  = \sum_{\substack{m,n\geq 1\\mn\leq x}} g(m) =\sum_{1\leq m\leq x} g(m)
  \sum_{1\leq n\leq x/m} \!1=
  \sum_{1\leq m\leq x} g(m) \big(\frac{x}{m}+\bigO(1)\big)\\ & =
  x\sum_{m=1}^\infty \frac{g(m)}{m} +
  \bigO\big(x\sum_{m\geq x} \frac{g(m)}{m}\big)+
  \bigO\big(\sum_{1\leq m\leq x} g(m)\big)=
  x\sum_{m=1}^\infty \frac{g(m)}{m}+\bigO(1)\\ & =
  x\prod_{p {\rm ~prime}}\big(1+\frac{1}{p^2(p^2-2)}\big)+\bigO(1)=
  C_1\;x+\bigO(1)\;.
\end{align*}
By summation by parts, we hence have
\begin{align*}
  \sum_{1\leq n\leq x}n^3f(n) &=\int_1^xt^3d[S(t)]=
  \big[t^3(C_1\;t+\bigO(1))\big]_1^x-
  3\int_1^xt^2(C_1\;t+\bigO(1))\;dt\\ &=\frac{C_1}{4}\;x^4+\bigO(x^3)\;.
\end{align*}
This proves Equation \eqref{eq:asympcube} and concludes the proof of
Proposition \ref{prop:Fouvry2}.  \cqfd

\section{Pair correlations of common perpendiculars in the
  modular curve $\PSLZ\bs\hdr$}
\label{sect:geometricmotivation}

In this section, we give a geometric motivation for the introduction
of the Euler function as multiplicities in the family $\L_\NN^\varphi$ of
logarithms of natural numbers. We refer to
\cite{ParPau17ETDS,BroParPau19} for more information.

Let $Y$ be a nonelementary geodesically complete connected proper
locally $\CAT(-1)$ good orbispace, so that the underlying space of $Y$
is $\Ga\bs\wt Y$ with $\wt Y$ a geodesically complete proper
$\CAT(-1)$ space and $\Ga$ a discrete group of isometries of $\wt Y$
preserving no point nor pair of points in $\wt Y\cap\partial_\infty
\wt Y$.  Let $D^-$ and $D^+$ be connected proper nonempty properly
immersed locally convex closed subsets of $Y$, that is, $D^-$ and
$D^+$ are the $\Ga$-orbits of proper nonempty closed convex subsets
$\wt D^-$ and $\wt D^+$ of $\wt Y$.  A {\it common perpendicular}
$\alpha$ between $D^-$ and $D^+$ is the $\Ga$-orbit of the unique
shortest arc $\wt \alpha$ between $\wt D^-$ and $\ga \wt D^+$ for some
$\ga\in\Ga$ such that $d(\wt D^-,\ga \wt D^+)>0$. The {\it multiplicity}
$\mult(\alpha)$ of $\alpha$ is the ratio $A/B$ where

$\bullet$~ $A$ is the number of elements $(\ga_-,\ga_+)\in
(\Ga/\Ga_{D^-}) \times (\Ga/\Ga_{\ga D^+}) $ such that $\wt \alpha$ is
the unique shortest arc between $\ga_-\wt D^-$ and $\ga_+\ga \wt D^+$,
and

$\bullet$~ $B$ is the cardinality of the pointwise stabilizer of
$\wt \alpha$ in $\Ga$.

\noindent The {\it length} $\lambda(\alpha)$ of the common
perpendicular $\alpha$ is the length of the geodesic segment $\wt
\alpha$ in $\wt Y$.  For every $\ell$ in the set
$\operatorname{OL}^\natural(D^-,D^+)$ of lengths of common
perpendiculars, the {\it length multiplicity} of $\ell$ is the sum of
the multiplicities of the common perpendiculars between $D^-$, $D^+$
having the length $\ell$ :
$$
\omega(\ell)=\sum_{\substack{\alpha{\rm~common~perpendicular}\\
     {\rm beween~} D^-{\rm~and~} D^+{\rm~with~}\lambda(\alpha)=\ell}}
\mult(\alpha)\;.
$$
If $\Perp(D^-,D^+)$ is the set of all common perpendiculars from $D^-$
to $D^+$ with multiplicities, then $(\lambda(\alpha))
_{\alpha\in\Perp(D^-,\,D^+)}$ is the {\em marked ortholength spectrum}
from $D^-$ to $D^+$, and the set $\ols(D^-,\,D^+)=
(\operatorname{OL}^\natural(D^-,D^+),\;\omega)$ of the lengths of the
common perpendiculars endowed with the length multiplicity $\omega$ is
the {\em ortholength spectrum} from $D^-$ to $D^+$.

The {\it pair correlation measure of the common perpendiculars from
$D^-$ to $D^+$} is the pair correlation measure of
$$\F_{D^-,D^+}=
\big((F_N=\operatorname{OL}^\natural(D^-,D^+)\cap[0,2\ln N])_{N\in\NN},\;
\omega\big)\,.
$$
We will study the asymptotic properties of the pair correlation
measures $\R^{\F_{D^-,D^+},\,\psi}_N$ for appropriate scalings $\psi$
in a subsequent paper, and we only consider in this paper the
following example.

\medskip
Let
$$
\wt Y=\hdr= \big(\{z\in\CC : \Im\;z>0\},ds^2=
\frac{d(\Re\;z)d(\Im\;z)}{(\Im\;z)^2}\big)
$$
be the upper halfspace model of the real hyperbolic plane with
constant curvature $-1$. For every $b\in\NN-\{0\}$, let $\Ga_0[b]$ be
{\it Hecke's congruence subgroup modulo $b$} of the {\it modular
  group} $\PSLZ$, which is the preimage of the upper triangular
subgroup $\PSL(\ZZ/b\ZZ)$ under the reduction morphism
$\PSLZ\ra\PSL(\ZZ/b\ZZ)$. It acts faithfully by homographies on
$\hdr$, and is a lattice in the isometry group of $\hdr$. Let
$Y^b=\Ga_0[b]\bs\hdr$, which is a finite (ramified) cover of the {\it
  modular curve} $\PSLZ\bs\hdr$.  Let $\wt D^-=\wt D^+$ be the
horoball $\H_\infty=\{z\in\CC : \Im\;z\geq 1\}$ in $\hdr$, whose image
$D^-=D^+$ in $Y^b$ is a {\it Margulis neighbourhood} of a cusp of
$Y^b$. If $b=1$, then $D^-=D^+$ is actually the maximal Margulis
neighbourhood of the unique cusp of $Y^b$. To emphasize the dependence
on the integer $b$, we will use the notation $\F^b_{D^-,D^+}=
\F_{D^-,D^+}$ for the family of lengths of common perpendiculars
between $D^-$ and $D^+$

The following result says that the pair correlation measures of the
common perpendiculars from this Margulis cusp neighbourhood to itself
are, up to the homothety of factor $2$, the pair correlation measures
of the logarithms of the natural numbers congruent to $0$ modulo $b$,
with multiplicities given by the Euler function $\varphi$.

We use in the following result the notation
$\R^{\L_\NN^{\varphi,b,b},\,\psi}_N$ of the introduction with
$$
\L_\NN^{\varphi,b,b}=(\,(L_N=\{\ln n\;:\; 0<n\leq N,
\;n\equiv 0\!\!\mod b\})_{N\in\NN},\; \omega=\varphi\circ \exp)\;.
$$

\bprop\label{prop:corrhoromodular}
For every scaling function $\psi$ and every $N\in\NN$, if $f:t\mapsto
2t$, then $\R^{\F^b_{D^-,D^+},\,\psi}_N=f_*\big(\R^{\L_\NN^{\varphi,b,b},\,\psi}_N\,\big)$.
\eprop

\dem The orbit of $\H_\infty$ under $\Ga_0[b]$ consists, besides
$\H_\infty$ itself, of the Euclidean disks $\H_{\frac{p}{q}}$ of
Euclidean radius $\frac 1{2q^2}$ tangent to the horizontal line at the
rational numbers $\frac pq$ with $q>0$, $q\equiv 0\!\!\mod b$ and
$(p,q)=1$.

\medskip
\noindent \begin{minipage}{8cm} ~~~ Every common perpendicular between
  $D^-$ and $D^+$ has a unique representative which starts from the
  Euclidean segment $[i,i+1[$ on the boundary of $\H_\infty$ and ends
  on the boundary of $\H_{\frac{p}{q}}$ with $\frac{p}{q}\in\QQ\cap[0,1[$
  and $q\equiv 0\!\!\mod b$.  Its hyperbolic length is $2\ln q$.  In
  particular, we have $\operatorname{OL}^\natural(D^-,D^+) =\{2\ln q :
  q\geq 2,\;q\equiv 0\!\!\mod b\}$.
\end{minipage}
\begin{minipage}{6.8cm}
\:\:\:\begin{picture}(0,0)%
\includegraphics{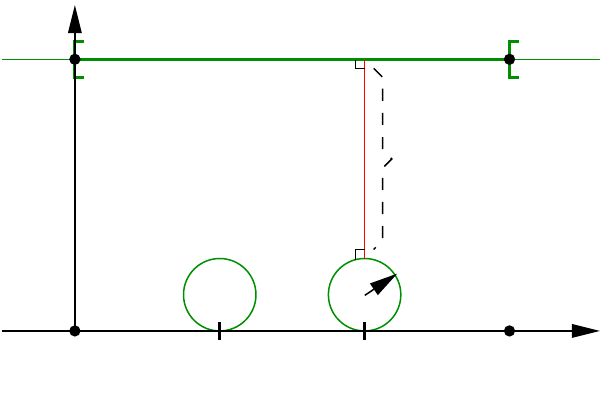}%
\end{picture}%
\setlength{\unitlength}{3812sp}%
\begingroup\makeatletter\ifx\SetFigFont\undefined%
\gdef\SetFigFont#1#2#3#4#5{%
  \reset@font\fontsize{#1}{#2pt}%
  \fontfamily{#3}\fontseries{#4}\fontshape{#5}%
  \selectfont}%
\fi\endgroup%
\begin{picture}(2994,1971)(1159,-1300)
\put(1486,-1141){\makebox(0,0)[lb]{\smash{{\SetFigFont{11}{13.2}{\rmdefault}{\mddefault}{\updefault}{\color[rgb]{0,0,0}$0$}%
}}}}
\put(2206,-1231){\makebox(0,0)[lb]{\smash{{\SetFigFont{11}{13.2}{\rmdefault}{\mddefault}{\updefault}{\color[rgb]{0,0,0}$\frac{1}{q}$}%
}}}}
\put(2926,-1231){\makebox(0,0)[lb]{\smash{{\SetFigFont{11}{13.2}{\rmdefault}{\mddefault}{\updefault}{\color[rgb]{0,0,0}$\frac{p}{q}$}%
}}}}
\put(3151,-106){\makebox(0,0)[lb]{\smash{{\SetFigFont{11}{13.2}{\rmdefault}{\mddefault}{\updefault}{\color[rgb]{0,0,0}$2\ln q$}%
}}}}
\put(3781,434){\makebox(0,0)[lb]{\smash{{\SetFigFont{11}{13.2}{\rmdefault}{\mddefault}{\updefault}{\color[rgb]{0,0,0}$i+1$}%
}}}}
\put(1396,434){\makebox(0,0)[lb]{\smash{{\SetFigFont{11}{13.2}{\rmdefault}{\mddefault}{\updefault}{\color[rgb]{0,0,0}$i$}%
}}}}
\put(3691,-1141){\makebox(0,0)[lb]{\smash{{\SetFigFont{11}{13.2}{\rmdefault}{\mddefault}{\updefault}{\color[rgb]{0,0,0}$1$}%
}}}}
\put(3151,-691){\makebox(0,0)[lb]{\smash{{\SetFigFont{11}{13.2}{\rmdefault}{\mddefault}{\updefault}{\color[rgb]{0,0,0}$\frac{1}{2q^2}$}%
}}}}
\put(2251,524){\makebox(0,0)[lb]{\smash{{\SetFigFont{11}{13.2}{\rmdefault}{\mddefault}{\updefault}{\color[rgb]{0,0,0}$\H_\infty$}%
}}}}
\end{picture}%

\end{minipage}

\medskip
Since $\PSL_2(\RR)$ acts simply transitively on the unit tangent
bundle of $\hdr$, the multiplicities of the common perpendiculars are
equal to $1$. Hence the length multiplicity of $2\ln q$ is exactly the
number of elements $p\in\ZZ/q\ZZ$ coprime with $q$, that is,
$\omega(2\ln q) =\varphi(q)$.
\cqfd

\bigskip
The following results, computing the pair correlation functions at
trivial or linear scaling of the lengths of the common perpendiculars
from the Margulis cusp neighbourhood at infinity to itself in Hecke's
modular curve $\Ga_0[b]\bs\hdr$, follow immediately from Theorems
\ref{theo:logpaircorrelphi} and \ref{theo:logpaircorrelphipsi} with
$a=b$, which also give an error term, using Proposition
\ref{prop:corrhoromodular}.

\bcoro \label{coro:comperphoromodular} (1) For all $b$, as
$N\ra+\infty$, the pair correlation measures
$\R^{\F^b_{D^-,D^+},\,1}_N$ on $\RR$, renormalized to be probability
measures, weak-star converge to a measure absolutely continuous with
respect to the Lebesgue measure on $\RR$, with pair correlation
function given by $s\mapsto \frac{1}{2}\;e^{-\,\,|s|}$.

\medskip\noindent (2) For all $b$, as $N\ra+\infty$, the pair correlation
measures $\frac{1}{N^3}\R^{\F^b_{D^-,D^+},\,\id}_N$ on $\RR$ weak-star
converge to a measure absolutely continuous with respect to the
Lebesgue measure on $\RR$, with pair correlation function given by
$s\mapsto \frac{1}{s^4}{\displaystyle \sum_{1\leq
    k\leq\lfloor\,|s|\,\rfloor,\; k\equiv 0\!\!\!\mod b}} c_{b,b,k}
\;k^3$, where $c_{b,b,k}$ is defined in Equation \eqref{eq:cabk}. \cqfd
\ecoro

%With the notation of the end of Section \ref{sect:weightlogpaircorrel},
%note that $s\mapsto \wt\D_N\big(\frac{s}{2}\big)$ is, by
%Proposition \ref{prop:corrhoromodular}, the cumulative distribution
%function at time $N$ of the differences of the lengths of the common
%perpendiculars from the maximal Margulis cusp neighbourhood to itself
%in the modular curve $\PSLZ\bs\hdr$, that is, for all $s,s'\in \RR$
%with $s<s'$, we have
%$$
%\frac{\R^{\F_{D^-,D^+},\,1}_N}
%{\big\|\R^{\F_{D^-,D^+},\,1}_N\big\|}\,(\,]s,s'\,]) =
%\wt\D_N\big(\frac{s'}{2}\big)-\wt\D_N\big(\frac{s}{2}\big)\;.
%$$
%See the picture at the end of Section \ref{sect:weightlogpaircorrel}
%for an illustration

%\newpage
\appendix
\section{Appendix : A Mirsky formula with congruences, \\ by Etienne Fouvry}
\label{appendixFouvry}

The aim of this appendix is to give a proof of a version with
congruences of an arithmetic formula due to Mirsky \cite{Mirsky49},
which is used in the proof of Theorem \ref{theo:logpaircorrelphipsi}.

Let $a,b,k\in\NN$ be fixed integers satisfying $a,b\geq 1$.  Denoting
by $\varphi$ the Euler function, we give an asymptotic formula, as
$x\geq 1$ tends to $+\infty$, for the quantity
$$
S(x;a,b,k)= \sum_{\substack{1\leq n\leq x\\n\,\equiv\, a\!\!\!\!\mod b}}
\varphi(n)\;\varphi(n+k)\;.
$$
Throughout this appendix, the letter $p$ denotes as usual a varying
positive prime in $\ZZ$, and we use the convention on $\bigO(\cdot)$
of the beginning of Lemma \ref{lem:Fouvry1}. For all integers
$\alpha,\beta\geq 1$, we denote as usual by $(\alpha,\beta)$ and
$[\alpha,\beta]$ their positive gcd et lcm, respectively.

\btheo\label{theo:mirskyfouvry}
Let
\begin{equation}\label{eq:cabk}
c_{a,b,k}= \frac{1}{b}\;\prod_{\substack{p\\(p,b)\,\mid\, k+a}}
\big(1-\frac{(p,b)}{p^2}\big)\;
\prod_{\substack{p \\ (p,\,b) \,\mid\, a}}
\big(1-(p,b)\,\frac{\kappa_{a,b,k}(p)\;\kappa'_k(p)}{p^2}\big)\;,
\end{equation}
where
$$
\kappa_{a,b,k}(p)=\left\{\begin{array}{l}
(1-\frac{(p,b)}{p^2})^{-1}{\rm ~~~if~~~} (p,b)\mid a+k\\ 1
{\rm ~~~otherwise} \end{array}\right. {\rm ~~~and~~~}
\kappa'_{k}(p)=\left\{\begin{array}{l}
1-\frac{1}{p}{\rm ~~~if~~~} p\mid k\\ 1 {\rm ~~~otherwise.}
\end{array}\right.
$$
There exists an absolute constant $K>0$ such that for all integers
$a,b\geq 1$ and $k\geq 0$ and real number $x\geq 1$, we
have
$$
\big|\;S(x;a,b,k)-
c_{a,b,k}\,\big(\,\frac{x^3}{3}+\frac{k\,x^2}{2}\big)\;\;\big|\leq
K\big(x(x+k)(\ln 2x)\ln (2x+k)\big)\;.
$$
\etheo

\brema\label{rem:majominocabk}{\rm 
Before proving  this theorem, we give some considerations on the constant
$c_{a,b,k}.
$
We start from the obvious inequalities, for every prime $p$,
$$
1 \leq \kappa_{a,b,k} (p) \leq 2 \text{ ~~and~~ }
\frac{1}{2} \leq { \kappa'}_k (p) \leq 1\;.
$$

%From these inequalities and the definition \eqref{eq:cabk}, and by
%a discussion on the cases $a,b,k$ even or not in order to obtain a
%minoration by $\frac{1}{4}$ of the local factor corresponding to $p=2$
%in the Euler products, we deduce that $c_{a,b,k}$ satisfies the
%following inequalities
%\begin{equation}\label{eq:cabkc}
%\frac{1}{4\,b}\, \prod_{p\geq 3 \atop (p,b) \mid a+k}
%\Big(1-\frac{(p,b)}{p^2}\Big)\, \prod_{p\geq 3\atop (p,b)\mid a}
%\Big(1-2\cdot \frac {(p,b)}{p^2}\Big) \leq c_{a,b,k} < \frac 1b\;.
%\end{equation}
From these inequalities and the definition \eqref{eq:cabk},  we deduce
that $c_{a,b,k}$ satisfies the following inequalities
\begin{equation}\label{eq:cabkc}
\frac 1b \,\Lambda (a,b,k)\, \prod_{p\geq 3 \atop (p,b) \mid a+k}
\Big(1-\frac{(p,b)}{p^2}\Big)\, \prod_{p\geq 3\atop (p,b)\mid a}
\Big(1-2\cdot \frac {(p,b)}{p^2}\Big) \leq c_{a,b,k} < \frac 1b\;,
\end{equation}
where the local factor $\Lambda (a,b,k)$, obtained by separating the
case $p=2$ in the Euler products, is defined by
$$
\Lambda (a,b,k) =
\begin{cases} 5/8 &\text{ if } 2\nmid b \text{ and } 2 \mid k, \\
 1/2 &\text{ if } 2\nmid b \text{ and } 2 \nmid k, \\
 1 &\text{ if } 2 \mid b,\; 2 \nmid a \text{ and } 2 \mid k,\\
 1/2 &\text{ if } 2 \mid b,\; 2 \nmid a\text{ and } 2 \nmid k,\\
 1/4 &\text{ if } 2 \mid b,\; 2 \mid a \text{ and } 2 \mid k,\\
 1/2   &\text{ if } 2 \mid b,\; 2 \mid a \text{ and } 2 \nmid k\;.\\
 \end{cases}
$$
We have the positive lower bound $\Lambda (a,b,k)\geq \frac{1}{4}$ in
all cases. This shows in particular that $c_{a,b,k} >0$.  A deeper
look also leads to the statement that the positive product
$b\,c_{a,b,k}$ can be arbitrarily small : it suffices to consider the
case where all the integers $a$, $b$ and $k$ are all divisible by the
$t$ smallest primes, and letting the integer $t$ tend to infinity.  }
\end{rema}

Theorem \ref{theo:mirskyfouvry} was already known when $a=b=1$, this
result is due to \cite[Thm.~9, Eq.~(30)]{Mirsky49}, with
\begin{equation}\label{eq:c11kmirsky}
c_{1,1,k}= \prod_{p}\big(1-\frac{2}{p^2}\big)
\prod_{p\,\mid\, k}\big(1+\frac{1}{p(p^2-2)}\big)\;.
\end{equation}
For this computation, see Remark \ref{rem:retrouvMirsky} below.

\medskip
\noindent{{\bf Proof of Theorem \ref{theo:mirskyfouvry}. }
For every $x\geq 1$, let us first prove that there exists a constant
$c_{a,b,k}\in\;]0,1]$ such that the sum
$$
\wt S(x)=\sum_{\substack{1\leq n\leq x\\n\,\equiv\, a\!\!\!\!\mod b}}
\frac{\varphi(n)}{n}\cdot\frac{\varphi(n+k)}{n+k}
$$
satifies the asymptotic formula, uniformly in $a, b\geq 1$, $k\geq 0$
and $x\geq 1$,
\begin{equation}\label{eq:asymwtS}
\wt S(x)=c_{a,b,k}\;x+\bigO((\ln 2x)\ln (2x+k))\;.
\end{equation}
Theorem \ref{theo:mirskyfouvry} follows classically, by applying Abel's
summation formula
$$
\sum_{1\leq n\leq x}
a_nf(n)=\Big(\sum_{1\leq n\leq x} a_n\Big)\;f(x) -
\int_{1}^x\Big(\sum_{1\leq n\leq t} a_n\Big)\;f'(t)\;dt
$$
to the numerical sequence $\big(a_n=\frac{\varphi(n)}{n}\cdot
\frac{\varphi(n+k)}{n+k}\;\delta_n\big)_{n\geq 1}$, where $\delta_n=1$
if $n\equiv a\!\!\!\mod b$ and $\delta_n=0$ otherwise, and to the
function $f:[1,+\infty[\;\ra\RR$ of class ${\rm C}^1$ defined by
$x\mapsto x(x+k)$. We indeed have
\begin{align*}
S(x;a,b,k)&=\wt S(x) x(x+k)-\int_{1}^x\wt S(t)(2t+k)\;dt
\\ & =\frac{c_{a,b,k}}{3}\;x^3+\frac{c_{a,b,k}\,k}{2}\;x^2+
(\frac{2}{3}+\frac{k}{2})c_{a,b,k}+
\bigO\big(x(x+k)(\ln 2x)\ln (2x+k)\big)\;,
\end{align*}
which gives the result since $c_{a,b,k}\leq 1$.

\medskip
Let us denote by ${\bf 1}$ the constant arithmetic function with value
$1$. The convolution equality $\varphi=\mu\star \id$ implies by
division that $\frac{\varphi}{\id}=\frac{\mu}{\id}\star 1$. Applying
twice this formula, we have
$$
\wt S(x)
%=\sum_{1\leq d\leq x}
%\frac{\mu(d)}{d}\;\sum_{\substack{1\leq n\leq x,\;d\,\mid\, n\\
%n\,\equiv\, a\!\!\!\!\mod b}}\frac{\varphi(n+k)}{n+k}
=\sum_{1\leq d\leq x}
\frac{\mu(d)}{d}\;\sum_{1\leq \delta\leq x+k}
\frac{\mu(\delta)}{\delta}\;\sum_{\substack{1\leq n\leq x\\
\;d\,\mid\, n,\;\delta\,\mid\,(n+k)\\
n\,\equiv\, a\!\!\!\!\mod b}}1\;.
$$
The system of three congruences
${\displaystyle
n\equiv\left\{\begin{array}{l} 0\!\!\mod d\\
-k\!\!\mod \delta\\ a\!\!\mod b\end{array}
\right.
}$ has a solution $n\leq x$ if and only if there exists an integer
$m\leq x/d$ such that $n=dm$ and
\begin{equation}\label{eq:systdeucongru}
\left\{\begin{array}{l} dm\equiv 
-k\!\!\!\mod \delta\\  dm\equiv a\!\!\!\mod b\;.\end{array}
\right.
\end{equation}
When $(d,\delta)\nmid k$ or when $(d,b)\nmid a$, no solution
exists. We hence have
$$
\wt S(x)
=\sum_{\substack{1\leq d\leq x,\;1\leq \delta\leq x+k\\
(d,\delta)\,\mid\, k,\;(d,b)\,\mid\, a}}
\frac{\mu(d)}{d}\; \frac{\mu(\delta)}{\delta} \;\card
\big\{m\leq x/d : \begin{array}{l} dm\equiv -k\!\!\!\mod \delta\\
dm\equiv a\!\!\!\mod b\end{array}\big\}\;.
$$
If $(d,\delta)\mid k$ and $(d,b)\mid a$, let us denote by
$\overline{\frac{d}{(d,\delta)}}$ the multiplicative inverse of the
integer $\frac{d}{(d,\delta)}$ modulo $\frac{\delta}{(d,\delta)}$ and by
$\overline{\frac{d}{(d,b)}}$ the multiplicative inverse of the integer
$\frac{d}{(d,b)}$ modulo $\frac{b}{(d,b)}$. The system of two
congruences \eqref{eq:systdeucongru}, after division of its first
equation by $(d,\delta)$ and of its second equation by $(d,b)$, is
then equivalent to the system
$$
\left\{\begin{array}{l} m\equiv
-\frac{k}{(d,\delta)}\;\overline{\frac{d}{(d,\delta)}}\!\!
\mod \frac{\delta}{(d,\delta)}\medskip\\
m\equiv \frac{a}{(d,b)}\;\overline{\frac{d}{(d,b)}}\!\!
\mod \frac{b}{(d,b)}\;.\end{array} \right.
$$
This system has a solution if and only if the following divisibility
condition holds
\begin{align*}
&\Big(\frac{\delta}{(d,\delta)},\frac{b}{(d,b)}\Big)\; \mid
\frac{k}{(d,\delta)}\;\overline{\frac{d}{(d,\delta)}}+
\frac{a}{(d,b)}\;\overline{\frac{d}{(d,b)}}\\
\Leftrightarrow \;\;&
\Big(\frac{\delta}{(d,\delta)},\frac{b}{(d,b)}\Big)\; \mid
\frac{k}{(d,\delta)}\;\frac{d}{(d,b)}+
\frac{a}{(d,b)}\;\frac{d}{(d,\delta)}\end{align*}
Since the integers $\frac{d}{(d,b)}$ and $\frac{d}{(d,\delta)}$ are
coprime with the gcd $\big(\frac{\delta}{(d,\delta)}
, \frac{b}{(d,b)}\big)$, we deduce by multiplication that the above
condition holds if and only if we have
$$
\big(\delta(d,b), b(d,\delta)\big)\; \mid
d(k+a)\;.
$$
The following  lemma, where $\bigO(1)$ is uniformly bounded in 
$\alpha,\beta,\alpha_0,\beta_0$, is elementary.

\blemm
For all integers $\alpha,\beta,\alpha_0,\beta_0\geq 1$ and real number
$y\geq 1$, we have
$$
\card \{m\leq y: m\equiv \alpha_0\!\!\!\mod \alpha{\rm ~and~}
m\equiv \beta_0\!\!\!\mod \beta\}=\left\{\begin{array}{l} 
0 {\rm ~if~} \alpha_0\nequiv\beta_0\!\!\!\mod (\alpha,\beta)\\
\frac{y}{[\alpha,\beta]}+\bigO(1) {\rm ~otherwise}\;. \;\;\;
\Box\end{array}\right.
$$
\elemm

This lemma implies that $\card\big\{m\leq x/d :
\text{\small $\begin{array}{l} dm\equiv -k\!\!\!\mod \delta\\
dm\equiv a\!\!\!\mod b\end{array}$}\big\}\;=
\frac{x}{d\big[\frac{\delta}{(d,\delta)}, \frac{b}{(d,b)}\big]}
\;+\bigO(1)$ under the assumption that $(d,\delta)\mid k$ and
$(d,b)\mid a$, where $\bigO(1)$ is uniformly bounded in $a,b,k,d,\delta$.
By the classical majoration of the harmonic series, we have
$$
\sum_{1\leq d\leq x,\;1\leq \delta\leq x+k} \frac{1}{d}\; \frac{1}{\delta}
=\bigO((\ln 2x)\ln(2x+k))\;.
$$
By the relation between lcm, gcd and product of two positive integers,
we hence have
$$
\wt S(x) =\;x
\sum_{\substack{1\leq d\leq x,\;1\leq \delta\leq x+k\\
(d,\delta)\,\mid\, k,\;(d,b)\,\mid\, a\\
(\delta(d,b), b(d,\delta))\, \mid\, d(k+a)}}
\frac{\mu(d)}{d}\; \frac{\mu(\delta)}{\delta}\;
\frac{(\delta(d,b),\, b(d,\delta))}{d\delta b}
\;+\bigO((\ln 2x)\ln(2x+k))\;, 
$$
uniformly in $a,b\geq 1$, $k\geq 0$ and $x\geq 1$.

Completing the sum with the indices $d> x$ and $\delta> x+k$
introduces an error of the form (uniformly in $a,b\geq 1$, $k\geq 0$
and $x\geq 1$)
$$
\bigO\Big(\sum_{d\geq x,\;\delta\geq 1}
\frac{(d,\delta)}{d^2\delta^2}\Big)
=\bigO\Big(\sum_{t\geq 1,\;d'\geq x/t,\;\delta'\geq 1}
\frac{t}{(td')^2(t\delta')^2}\Big)=\bigO\big(\frac{1}{x}\big)\;.
$$
This proves Formula \eqref{eq:asymwtS} by setting
\begin{equation}\label{eq:cabktech}
c_{a,b,k}=\sum_{\substack{d,\;\delta\geq 1\\
(d,\delta)\,\mid\, k,\;(d,b)\,\mid\, a\\
(\delta(d,b), b(d,\delta))\, \mid\, d(k+a)}}
\frac{\mu(d)}{d}\; \frac{\mu(\delta)}{\delta}\;
\frac{(\delta(d,b),\, b(d,\delta))}{d\delta b}\;. 
\end{equation}

Let us now prove Equation \eqref{eq:cabk}. By
Remark \ref{rem:majominocabk}, this implies that $0<c_{a,b,k}\leq 1$,
hence completes the proof of Theorem \ref{theo:mirskyfouvry}.

\bigskip
\noindent{{\bf Proof of Equation \eqref{eq:cabk}. }
For every integer $d\geq 1$, let $\chi_d$ be the characteristic
function of the set of integers $\delta\geq 1$ such that $(\delta,d)\mid k$.
For every integer $d\geq 1$, let us define
\begin{equation}\label{eq:defpsisubd}
\psi_d: \delta\mapsto \big(\delta, \frac{b}{(d,b)}\,(d,\delta)\big)\;.
\end{equation}
Note that the assertion $(\delta(d,b), b(d,\delta)) \mid d(k+a)$
is equivalent to the assertion
$$
\psi_d(\delta)\, \mid\, \frac{d}{(d,b)}\,(k+a)\;.
$$
For every integer $d\geq 1$, let $\chi^*_d$ be the characteristic
function of the set of integers $\delta\geq 1$ such that the above
divisibility assertion is satisfied.  Let us define
\begin{equation}\label{eq:defceuler}
c^*:d\mapsto\sum_{\delta\geq 1}\frac{\mu(\delta)}{\delta^2}
\;\chi_d(\delta)\;\chi^*_d(\delta)\;\psi_d(\delta)
\end{equation}
(this arithmetic function $c^*$ depends on the constants
$a,b,k$). Equation \eqref{eq:cabktech} then becomes
\begin{equation}\label{eq:cabktech2}
c_{a,b,k}=\;\frac{1}{b}\sum_{\substack{d\geq 1\\(d,b)\,\mid\, a}}
\frac{\mu(d)}{d^2}\; (d,b)\;c^*(d)\;.
\end{equation}
In order to transform the series $c^*(d)$ defined by
Formula \eqref{eq:defceuler} into an Eulerian product and in order to
analyse it, we will use the following two lemmas.

\blemm\label{lem:multiplicatives}
For every integer $d\geq 1$, the arithmetic functions  $\chi_d$,
$\chi^*_d$ and $\psi_d$ are multiplicative.
\elemm

\dem
We have $\chi_d(1)=\chi^*_d(1)=\psi_d(1)=1$. Let $\delta_1,\delta_2$
be two coprime integers.

The equality $(\delta_1\delta_2,d)=(\delta_1,d)(\delta_2,d)$ and the
fact that $(\delta_1,d)$ and $(\delta_2,d)$ are coprime imply the
multiplicativity of $\chi_d$.

In order to prove the multiplicativity of the function $\psi_d$, we
write
$$
\psi_d(\delta_1\delta_2) =
\big(\delta_1\delta_2, \frac{b}{(d,b)}\,(d,\delta_1\delta_2)\big)=
\big(\delta_1, \frac{b}{(d,b)}\,(\delta_1,d)(\delta_2,d)\big)
\big(\delta_2, \frac{b}{(d,b)}\,(\delta_1,d)(\delta_2,d)\big)\;.
$$
Since $\delta_1$ is coprime to $(\delta_2,d)$ and since $\delta_2$
is coprime to  $(\delta_1,d)$, we obtain as wanted the equality
$\psi_d(\delta_1\delta_2) = \psi_d(\delta_1)\,\psi_d(\delta_2)$.

Finally, the multiplicativity of the function $\chi^*_d$ is a
consequence of the multiplicativity of the function $\psi_d$ and of
the fact that $\psi_d(\delta_1)$ and $\psi_d(\delta_2)$ are coprime.
\cqfd

\blemm\label{lem:valeurschichippsi}
For every prime $p$ and every integer $d\geq 1$, we have
$$
\psi_d(p)= \left\{\begin{array}{l} p {\rm ~~~ if~~~} p\mid d,\\
(p,b){\rm ~~~ otherwise},\end{array}\right.
$$
and
$$
\chi_d(p)\;\chi^*_d(p)=1\Leftrightarrow \left\{\begin{array}{l}
p\mid(d,k) {\rm ~~~ and~~~} p\mid\frac{d}{(d,b)}\,(k+a),\\
{\rm or}\\
p \nmid d  {\rm ~~~ and~~~} (p,b)\mid k+a\;.\end{array}\right.
$$
\elemm

\dem
The first formula follows from the definition of $\psi_d(p)$ (see
Formula \eqref{eq:defpsisubd}) by considering the three cases ($p \mid
d$), ($p \nmid d$ and $p\mid b$), and ($p \nmid d$ and $p\nmid b$).

The second formula follows from the first one, from the definitions of
$\chi_d(p)$ and $\chi^*_d(p)$, and from the fact that
$\chi_d(p)\;\chi^*_d(p) =1$ if and only if $\chi_d(p)=\chi^*_d(p)=1$,
by considering the two cases ($p \mid d$) and ($p \nmid d$).
\cqfd

\medskip
The arithmetic function $\delta\mapsto\mu(\delta)\chi_d(\delta)
\;\chi^*_d(\delta) \;\psi_d(\delta)$ being multiplicative by
Lemma \ref{lem:multiplicatives}, and vanishing on the nontrivial
powers of primes, the series defining $c^*(d)$ in
Formula \eqref{eq:defceuler} may be written as an Eulerian product
\begin{equation}\label{eq:prodeuler}
c^*(d)=\prod_p\big(1-\frac{\chi_d(p)
\;\chi^*_d(p)\;\psi_d(p)}{p^2}\big)=\prod_{\substack{p\\\chi_d(p)
\;\chi^*_d(p)=1}}\big(1-\frac{\psi_d(p)}{p^2}\big)\;.
\end{equation}
By Equations \eqref{eq:cabktech2} and \eqref{eq:prodeuler}, and by
Lemma \ref{lem:valeurschichippsi}, we have
$$
c_{a,b,k}=\;\frac{1}{b}\sum_{\substack{d\geq 1\\(d,b)\,\mid\, a}}
\frac{\mu(d)}{d^2}\; (d,b)\;\prod_{\substack{p\,\nmid\,d\\(p,b)\,\mid\,k+a}}
\big(1-\frac{(p,b)}{p^2}\big)
\prod_{\substack{p\,\mid\, (d,k)\\p\,\mid\, \frac{d}{(d,b)}(k+a)}}
\big(1-\frac{1}{p}\big)\;.
$$
Let us define $\Ga_{a,b,k}={\displaystyle\prod_{\substack{p\\(p,b)\,\mid\, k+a}}}
\big(1-\frac{(p,b)}{p^2}\big)$, so that
\begin{equation}\label{eq:cabktech3}
c_{a,b,k}=\;\frac{\Ga_{a,b,k}}{b}\sum_{\substack{d\geq 1\\(d,b)\,\mid\, a}}
\frac{\mu(d)}{d^2}\; (d,b)\;\prod_{\substack{p\,\mid\, d\\(p,b)\,\mid\, k+a
}}\big(1-\frac{(p,b)}{p^2}\big)^{-1}
\prod_{\substack{p\,\mid\, (d,k)\\p\,\mid\, \frac{d}{(d,b)}(k+a)}}
\big(1-\frac{1}{p}\big)\;.
\end{equation}
For every integer $d\geq 1$ without square factor such that $(d,b)\mid
a$, we have
\begin{align*}
\prod_{\substack{p\,\mid\, (d,k)\\p\,\mid\, \frac{d}{(d,b)}(k+a)}}
\big(1-\frac{1}{p}\big)& =
\prod_{\substack{p\,\mid\, (d,k)\\p\,\mid\, \frac{d}{(d,b)}}}
\big(1-\frac{1}{p}\big)
\prod_{\substack{p\,\mid\, (d,k)\\p\,\mid\, k+a}}\big(1-\frac{1}{p}\big)
\prod_{\substack{p\,\mid\, (d,k)\\p\,\mid\, (\frac{d}{(d,b)},k+a)}}
\big(1-\frac{1}{p}\big)^{-1} \\ & =
\prod_{p\,\mid\, (\frac{d}{(d,b)},k)}
\big(1-\frac{1}{p}\big)
\prod_{p\,\mid\, (d,a,k)}\big(1-\frac{1}{p}\big)
\prod_{p\,\mid\, (\frac{d}{(d,b)},a,k)}
\big(1-\frac{1}{p}\big)^{-1}\\ & =
\prod_{p\,\mid\, (\frac{d}{(d,b)},k)}
\big(1-\frac{1}{p}\big)
\prod_{p\,\mid\, (d,a,b,k)}\big(1-\frac{1}{p}\big)\\& =
\prod_{p\,\mid\, (\frac{d}{(d,b)},k)}
\big(1-\frac{1}{p}\big)
\prod_{p\,\mid\, (d,b,k)}\big(1-\frac{1}{p}\big)=
\prod_{p\,\mid\, (d,k)}\big(1-\frac{1}{p}\big)\;.
\end{align*}

Thus, Equation \eqref{eq:cabktech3} writes $c_{a,b,k}$ as a series
$\frac{\Ga_{a,b,k}}{b}{\displaystyle \sum_{\substack{d\geq 1\\
(d,b)\,\mid\,a}}} \frac{f(d)}{d^2}$ where $f$ is a multiplicative
function, which vanishes on the nontrivial powers of primes. By
Eulerian product, we have therefore proved
Equation \eqref{eq:cabk}. \cqfd

\medskip
\brema\label{rem:retrouvMirsky}
{\rm When $a=b=1$, we indeed recover Mirsky's result \cite[Thm.~9,
Eq.~(30)]{Mirsky49}. Indeed, by Equation \eqref{eq:cabk}, we have
\begin{align*}
c_{1,1,k}&=\prod_{p}\big(1-\frac{1}{p^2}\big) \prod_{p\,\mid\, k}
\Big(1-\frac{(1-\frac{1}{p^2})^{-1}(1-\frac{1}{p})}{p^2}\Big)
\prod_{p\,\nmid\, k}\Big(1-\frac{(1-\frac{1}{p^2})^{-1}}{p^2}\Big)\\ & =
\prod_{p}\big(1-\frac{1}{p^2}\big)
\prod_{p\,\mid\, k}\Big(1-\frac{p-1}{p(p^2-1)}\Big)
\prod_{p}\big(1-\frac{1}{p^2-1}\big)
\prod_{p\,\mid\, k}\Big(1-\frac{1}{p^2-1}\Big)^{-1}\\ & =
\prod_{p}\big(1-\frac{2}{p^2}\big)
\prod_{p\,\mid\, k}\big(1+\frac{1}{p(p^2-2)}\big)\;.
\end{align*}}
\erema

%\newpage

{\small \bibliography{../biblio} }
%{\small \bibliography{/users/jouniparkkonen/Documents/CloudDocs/Documents/Latex/viitteet} }

\bigskip
{\small
\noindent \begin{tabular}{l} 
Department of Mathematics and Statistics, P.O. Box 35\\ 
40014 University of Jyv\"askyl\"a, FINLAND.\\
{\it e-mail: jouni.t.parkkonen@jyu.fi}
\end{tabular}
\medskip

\noindent \begin{tabular}{l}
Laboratoire de mathématique d'Orsay, UMR 8628 CNRS,\\
Universit\'e Paris-Saclay,\\
91405 ORSAY Cedex, FRANCE\\
{\it e-mail: frederic.paulin@universite-paris-saclay.fr}
\end{tabular}
\medskip

\noindent \begin{tabular}{l}
Laboratoire de mathématique d'Orsay, UMR 8628 CNRS,\\
Universit\'e Paris-Saclay,\\
91405 ORSAY Cedex, FRANCE\\
{\it e-mail: etienne.fouvry@universite-paris-saclay.fr}
\end{tabular}
}
\end{document}